\documentclass[final,onefignum,oneeqnum,onetabnum]{siamart171218}

 \usepackage{amsmath}
 \usepackage{float}
\usepackage{amsfonts}
\usepackage{amssymb}

\newcommand{\cS}{{\mathfrak S}}

\newcommand{\cD}{{\mathfrak D}}

 \newtheorem{remark}{Remark}[section]

\newcommand{\var}{{\rm Var}}

\newcommand{\Q}{{\cal Q}}

\newcommand{\A}{{\cal A}}

\newcommand{\be}{\begin{equation}}
\newcommand{\ee}{\end{equation}}

\newcommand{\ep}{{\epsilon}}

\def\e{\varepsilon}

\def\vv{\vartheta}

\def\bbr{{\Bbb{R}}} 
\def\bbe{{\Bbb{E}}} 

 \renewcommand{\theequation}{\thesection.\arabic{equation}}

\usepackage{amsfonts}
\usepackage{graphicx}
\usepackage{epstopdf}
\usepackage{algorithmic}
\usepackage{booktabs}
\ifpdf
  \DeclareGraphicsExtensions{.eps,.pdf,.png,.jpg}
\else
  \DeclareGraphicsExtensions{.eps}
\fi

\newcommand{\TheTitle}{Duality and sensitivity analysis of multistage linear stochastic programs}
\newcommand{\TheAuthors}{}

\author{
  Vincent Guigues\thanks{School  of Applied Mathematics, Funda\c{c}\~ao Getulio Vargas,
190 Praia de Botafogo, Rio de Janeiro, Brazil,
    (\email{vincent.guigues@fgv.br}). Research of this author was partially supported
    by CNPq grants 04872/2018-9 and 311289/2016-9.}
  \and
     Alexander Shapiro\thanks{
     School of Industrial and Systems Engineering,
Georgia Institute of Technology, Atlanta, GA 30332-0205, USA,
    (\email{ashapiro@isye.gatech.edu}).
    Research of this author  was partly supported by NSF grant 1633196.}
    \and
   Yi  Cheng \thanks{
     School of Industrial and Systems Engineering,
Georgia Institute of Technology, Atlanta, GA 30332-0205, USA,
    (\email{cheng.yi@gatech.edu}).}
}

\title{{\TheTitle}}

\usepackage{amsopn}

\ifpdf
\hypersetup{
  pdftitle={\TheTitle},
  pdfauthor={\TheAuthors}
}
\fi

\begin{document}




\maketitle


\begin{abstract} In this paper we investigate  the dual of a Multistage
Stochastic Linear Program (MSLP) to study two related questions for this class of problems. The first of these questions is the study of the
optimal value of the problem as a function of the involved  parameters.
For this sensitivity analysis problem, we provide formulas for the derivatives of the value function with respect to the parameters and illustrate their application
on an inventory problem. Since these formulas involve optimal dual solutions,
we need an algorithm that computes such solutions to use them, i.e., we need to solve
the dual problem.

In this context, the second question we address is the study of solution
methods for the dual problem. Writing Dynamic Programming equations for the dual, we
can use an SDDP type method, called Dual SDDP, which solves these Dynamic Programming equations
computing a sequence of nonincreasing deterministic upper bounds on the optimal value
of the problem. However, applying
this method will only be possible if the Relatively Complete Recourse (RCR) holds for the dual.
Since the RCR assumption may fail to hold (even for simple problems),
we design two variants of Dual SDDP, namely Dual SDDP with penalizations and Dual SDDP with
feasibility cuts, that converge to the optimal value of the dual (and therefore primal when
there is no duality gap) problem under mild assumptions. We also show that optimal dual
solutions can be obtained computing dual solutions of the subproblems solved when applying
Primal SDDP to the original primal MSLP.

The study of this second question allows us to take a fresh look at the
class of MSLP with interstage dependent cost coefficients. Indeed, for this
class of problems, cost-to-go functions are non-convex and solution methods were so far
using SDDP for a Markov chain approximation of the cost coefficients process.
For these problems, we propose to apply Dual SDDP with penalizations to the
cost-to-go functions of the dual which are concave.
This algorithm converges to the optimal value of the problem.

Finally, as a proof of concept of the tools developed, we present the results
of numerical experiments computing the sensitivity of the optimal value of an inventory
problem as a function of parameters of the demand process
and compare Primal and Dual SDDP on the inventory and a hydro-thermal planning problems.
\end{abstract}



\section{Introduction}\label{sec-intr}

Duality plays a key role in optimization. For generic optimization problems,
weak duality allows  to bound the optimal value.
Dual information is also used in many optimization algorithms such as
Uzawa algorithm \cite{uzawabook}, primal-dual projected gradient \cite{nedic} or Stochastic Dual Dynamic Programming (SDDP) \cite{pereira}. Moreover, for several classes of optimization problems, the dual is easier to solve than the primal problem, for instance
when it is amenable to decomposition techniques such as price decomposition \cite{bertsekas}.
Even when there is a duality gap between the primal and dual optimal values, solving the dual
already gives a bound on the optimal value, as mentioned earlier.
Duality is also a fundamental tool in the reformulation of Robust Optimization problems, see for instance \cite{nem2}.
Finally, derivatives of the value function of classes of optimization problems can be
related to optimal dual solutions, see \cite{Bonnans2000}, \cite{Rockafellar}
and more recently \cite{guiguessiopt2016, guiguesinexactsddp, guiguesinexactsmd} for
the characterization of subdifferentials, subgradients, and
$\varepsilon$-subgradients of value functions of convex optimization problems.

For stochastic control problems, stochastic Lagrange multipliers
were already used in \cite{kushner69, kushner65, kushner72}.
In the context of multistage stochastic programs, duality was studied in
\cite{rockduality99, higlesen2006}, see also \cite{shadenrbook} for a review.
More recently, the sensitivity analysis of multistage stochastic programs was discussed in \cite{bonnan}  and \cite{woza}.
In \cite{bonnan} the authors study the sensitivity with respect to parameters driving the considered  price model. The corresponding parameters are in the objective function and the analysis of  the estimate of marginal price is  based on  Danskin's theorem with the SDDP
method used for the numerical calculations.
In  \cite{woza}, the authors use the Envelope Theorem for the sensitivity analysis. The required derivatives are described in terms of Lagrange multipliers associated with the value functions.

In this paper, focusing our attention on the dual of a Multistage Stochastic Linear Program (MSLP), we are able to
provide insights into  three important problems for MSLPs: sensitivity analysis,
computation of a sequence of deterministic upper bounds on the optimal value which converges to the optimal value, and use of duality
to solve Dynamic Programming (DP) equations on the dual which are simpler to solve (in the sense that they have convex
cost-to-go functions) than primal DP equations for problems with interstage dependent cost coefficients. Our main contributions are summarized below.

\par {\textbf{Sensitivity analysis of MSLPs.}} We explain how to compute derivatives
of the optimal value, seen as a function of the problem parameters, of a
MSLP in terms of dual optimal solutions. Therefore, the construction of the dual
problem is essential for our approach, contrary to \cite{bonnan}.
With respect to the sensitivity analysis \cite{woza}, in our approach, we do not use value functions directly, which are not known and can only be approximated, but rather construct the dual problem which is solved by an SDDP type algorithm, called
Dual SDDP.

\par {\textbf{Writing Dynamic Programming equations for the dual problem.}} A simple but crucial
ingredient for our developments and subsequent analysis of solution methods for the dual problem of a MSLP
is to write DP equations for that dual problem. We are not aware of another paper with these equations.
However, a similar study was done in \cite{cermics2018}. More precisely, for a stochastic linear
control problem with uncertainty in the right-hand-side, in \cite{cermics2018}, DP equations are
written for the conjugate of the cost-to-go functions and using an SDDP type method for these DP equations,
a sequence of upper bounds on the MSP optimal value is constructed which is the sequence of conjugate of the approximate
first stage cost-to-go functions evaluated at the initial state $x_0$.
Our approach has the advantage of being much simpler: contrary to derivations in \cite{cermics2018} which
require some algebra, our DP equations can be immediately obtained from the dual problem formulation,
this latter being known (given in \cite{shadenrbook} for instance).
On top of that, we relax two assumptions made in \cite{cermics2018}: (a) the relatively complete recourse
assumption of the dual and (b) randomness in the right-hand-side of the constraints only and interstage
independent. The next three paragraphs describe how the scope of (a) and (b) was extended in our
analysis.

\par {\textbf{Dual SDDP for dual problems without relatively complete recourse.}}
In \cite{cermics2018}, it is assumed that the dual problem of the considered MSLP satisfies
an assumption (Assumption 3) stronger than relatively complete recourse. This assumption may not be easy
to check or may not be satisfied (for instance it is not satisfied for the inventory and hydro-thermal
problems considered in Section \ref{numer}). Therefore, it is desirable to extend the scope of Dual
SDDP in such a way that it can still compute a deterministic converging sequence of upper bounds without
this assumption. We present two variants of Dual SDDP that can do that: Dual SDDP with penalizations
and Dual SDDP with feasibility cuts.

\par {\textbf{Dual SDDP for dual problems with all problem data random.}} Our DP equations are
written for problems with uncertainty in all parameters. We explain how to apply Dual SDDP
for such problems that do not satisfy (b) above.

\par {\textbf{Dual SDDP for problems with interstage dependent cost coefficients.}} Finally, we also
relax assumption (b) considering problems having interstage dependent cost coefficients.
Writing DP equations for the corresponding dual problem, we can apply Dual SDDP algorithm to solve
these equations, which, interestingly, have concave cost-to-go functions whereas primal cost-to-go
functions are not convex. This is in sharp contrast with the solution methods proposed so far such as
\cite{bonnan, shaplondorf} which apply SDDP on the primal cost-to-go functions using a Markov
chain approximation of the cost coefficients process.\\

The outline of the paper is the following.
Our building blocks are elaborated in
Section \ref{sec-dual} where we write DP equations for the dual, we explain how to
build upper bounding functions for the cost-to-go functions of the dual using penalizations, and study the dynamics of Lagrange multipliers. Sensitivity analysis of MSLPs is conducted in Section  \ref{sec-sens}
while Dual SDDP and its variants are studied in Section \ref{sec:dualsddp}.
Finally, the results of numerical simulations testing the tools developed on an inventory and an
hydro-thermal problem are presented in Section \ref{numer}. The interested reader can find
and test the code of all implementations and of Primal and Dual SDDP for MSLPs
at {\url{https://github.com/vguigues/Dual_SDDP_Library_Matlab}} and {\url{https://github.com/vguigues/Primal_SDDP_Library_Matlab}}. Proofs are collected in the Appendix.

\setcounter{equation}{0}
\section{Duality of multistage linear stochastic  programs}\label{sec-dual}

\subsection{Writing Dynamic Programming equations for the dual}\label{wrdpequ}

Consider the  multistage linear stochastic program
\begin{equation}\label{eq-1}
  \begin{array}{cll}
    \min\limits
    _{x_t\ge 0} & \bbe\left[  \sum_{t=1}^T c_t^\top x_t\right]\\
    {\rm s.t.}& A_1 x_1=b_1,\\
    & B_t x_{t-1}+A_t x_t=b_t,\;t=2,...,T.
   \end{array}
\end{equation}
Here     vectors $c_t=c_t(\xi_t)\in \bbr^{n_t}$, $b_t =b_t(\xi_t)\in \bbr^{m_t}$  and matrices $B_t=B_t(\xi_t)$, $A_t=A_t(\xi_t)$ are functions of  random process $\xi_t\in \bbr^{d_t}$, $t=1,...,T$   (with $\xi_1$ being deterministic). We denote by $\xi_{[t]}=(\xi_1,...,\xi_t)$ the history of the data process up to time $t$ and by $\bbe_{|\xi_{[t]}}$ the corresponding conditional expectation.
The optimization in \eqref{eq-1} is performed over functions (policies)  $x_t=x_t(\xi_{[t]})$, $t=1,...,T,$  of the data process  satisfying the feasibility constraints.

The Lagrangian of problem \eqref{eq-1} is
\begin{equation}\label{lagr}
\begin{array}{l}
 L(x,\pi)=\bbe\left[  \sum_{t=1}^T c_t^\top x_t+\pi_t^\top (b_t-B_t x_{t-1}-A_t x_t)\right]
 \end{array}
\end{equation}
in variables\footnote{Note that since $\xi_1$ is deterministic, the first stage decision  $x_1$  is  also deterministic; we write it as $x_1(\xi_{[1]})$ for uniformity of notation, and similarly for $\pi_1$.}
$x=(x_1(\xi_{[1]}),\ldots,x_T(\xi_{[T]}))$ and
$\pi=(\pi_1(\xi_{[1]}),\ldots,\pi_{T}(\xi_{[T]}))$
with the convention that $x_0=0$.
Dualization of the feasibility constraints leads to the following dual of problem \eqref{eq-1} (cf., \cite[Section 3.2.3]{shadenrbook}):
\begin{equation}\label{eq-2}
\begin{array}{cll}
  \max\limits_{\pi} & \bbe\big[  \sum_{t=1}^T     b_t^\top \pi_t\big]\\
    {\rm s.t.}  &A^\top_T\pi_T\le c_T,\\
  &  A_{t -1}^\top \pi_{t-1}+
    \bbe_{|\xi_{[t-1]}}\left[ B_{t}^\top \pi_{t}\right ]\le c_{t-1},\;t=2,...,T.
   \end{array}
\end{equation}
The optimization in \eqref{eq-2} is over policies $\pi_t=\pi_t(\xi_{[t]})$, $t=1,...,T$.

Unless stated otherwise, we make the  following assumption throughout the paper.
\begin{itemize}
   \item[(A1)]
    The process $\xi_1,...,\xi_T$ is stagewise independent  (i.e., random vector  $\xi_{t+1}$ is independent of $\xi_{[t]}$, $t=1,...,T-1$),
    and distribution of $\xi_t$ has a finite support, $\{\xi_{t 1},\ldots,\xi_{t N_t}\}$ with respective probabilities $p_{tj}$,
    $j=1,...,N_t$,  $t=2,...,T$.
 We denote by $A_{t j},B_{t j},c_{t j}, b_{t j}$ the respective  scenarios corresponding to $\xi_{t j}$.
 \end{itemize}

 Since  the random process $\xi_t$, $t=1,...,T$, has a {\em finite} number of realizations (scenarios),  problem \eqref{eq-1} can be viewed as a large linear program  and \eqref{eq-2} as its dual. By the standard theory of linear programming we have the following.

 \begin{proposition}
 \label{pr-duality}
 Suppose   that problem \eqref{eq-1} has a finite optimal value. Then   the optimal values of problems \eqref{eq-1} and \eqref{eq-2} are equal to each other and both problems have optimal solutions.
 \end{proposition}

We can write the following dynamic programming equations for the dual problem \eqref{eq-2}.  At the last stage $t=T$,
given $\pi_{T-1}$ and
$\xi_{[T-1]}$, we need to solve the following  problem with respect to $\pi_T$:
\begin{equation}\label{eq-3-1}
\begin{array}{cll}
  \max\limits_{\pi_T} &    \bbe[  b_T^\top \pi_T]\\
    {\rm s.t.}&
  A_{T}^\top \pi_{T}
    \le c_{T},\\
    &A_{T -1}^\top \pi_{T-1}+
    \bbe \left[ B_{T}^\top \pi_{T}\right ]\le c_{T-1}.
   \end{array}
\end{equation}
Since $\xi_T$  is independent of $\xi_{[T-1]}$,  the expectation in \eqref{eq-3-1} is unconditional with respect to the distribution of $\xi_T$.
In terms of scenarios the above problem can be written as
\begin{equation}\label{eq-3a-1}
\begin{array}{cll}
  \max\limits_{\pi_{T 1},\ldots,\pi_{T N_T}} &  \sum\limits_{j=1}^{N_T}   p_{T j}  b_{T j}^\top \pi_{T j}\\
    {\rm s.t.}&
  A_{T j}^\top \pi_{T j}  \le c_{T j},\;j=1,...,N_T,\\
    &A_{T -1}^\top \pi_{T-1}+
  \sum\limits_{j=1}^{N_T}   p_{T j}  B_{T j}^\top
     \pi_{T j} \le c_{T-1}.
   \end{array}
\end{equation}

The optimal value $V_T(\pi_{T-1},\xi_{T-1})$ and an optimal   solution\footnote{Note that problem \eqref{eq-3a-1} may have more than one optimal solution. In case of finite number of scenarios the considered linear program always has a solution provided its optimal value is finite.}
$(\bar{\pi}_{T 1},\ldots,\bar{\pi}_{T N_T})$  of   problem \eqref{eq-3a-1}  are  functions of vectors  $\pi_{T-1}$ and $c_{T-1}$ and matrix  $A_{T-1}$.
And  so on going backward in time, using  the stagewise independence assumption,  we can write the respective dynamic programming equations for
  $t=T-1,...,2$, as
 \begin{equation}\label{eq-6-1}
\begin{array}{cll}
  \max\limits_{\pi_{t 1},\ldots,\pi_{t N_t}} &     \sum\limits_{j=1}^{N_{t}}  p_{t j}  \left[ b_{t j}^\top \pi_{t j}+ V_{t+1}(\pi_{t j},\xi_{t j})\right]\\
    {\rm s.t.}&
  A_{t-1}^\top \pi_{t-1}+
  \sum\limits_{j=1}^{N_{t}} p_{tj} B_{t j}^\top  \pi_{t j}  \le c_{t-1},
   \end{array}
\end{equation}
with $V_t(\pi_{t-1},\xi_{t-1})$ being the optimal value  of problem \eqref{eq-6-1}.
Finally at the first stage the following problem should be solved
 \begin{equation}\label{eq-7-1}
  \max_{\pi_{1}}      b_1^\top \pi_1+ V_{2}(\pi_1,\xi_1) .
   \end{equation}

These dynamic  programming equations can be compared with the
dynamic  programming equations for primal problem \eqref{eq-1}, where the respective cost-to-go (value) function  $Q_t(x_{t-1},\xi_{t j})$, $j=1,...,N_t$,  is  given by the optimal value of
\begin{equation}\label{primal}
\begin{array}{cll}
\min\limits_{x_t\ge 0} & c_{t j}^\top x_t+\Q_{t+1}(x_t)\\
{\rm s.t.} & B_{t j} x_{t-1}+A_{t j} x_t=b_{t j},
\end{array}
\end{equation}
with
\[
\Q_{t+1}(x_t):=\bbe[ Q_{t+1}(x_t,\xi_{t+1})]=\sum_{j=1}^{N_{t}} p_{t+1 j} Q_{t+1}(x_t,\xi_{t+1 j}).
\]

Let us make the following observations  about   the dual problem.
\begin{itemize}
  \item [(i)]
  Unlike in  the primal problem,
  the optimization (maximization) problems \eqref{eq-3a-1} and \eqref{eq-6-1}   do not decompose into separate problems with respect to each $\pi_{t j}$ and should be solved as one linear program with respect to $(\pi_{t 1},...,\pi_{t N_t})$.
  \item [(ii)]
 The value function   $V_t(\pi_{t-1},\xi_{t-1})$   is a concave function of $\pi_{t-1}$.
 \item [(iii)]
If $A_t$ and $c_t$, $t=2,...,T$,  are deterministic, then $V_t(\pi_{t-1})$ is only a function  of $\pi_{t-1}$.
\end{itemize}

\subsection{Relatively complete recourse}\label{sec:rcr}

The following definition of Relatively Complete Recourse (RCR) is applied to the dual problem.    Recall that we assume that the set of possible realizations (scenarios) of the data process is finite.
\begin{definition}
\label{def-rcr}
We say that a sequence $\bar{\pi}_t$, $t=1,...,T$, is generated by the forward (dual)  process if $\bar{\pi}_1\in \bbr^{m_1}$ and   for    $\pi_{t-1}=\bar{\pi}_{t-1}$,   $t=2,...,T$, going forward in time,  $\bar{\pi}_t$  coincides with some $\pi_{t j}$, $j=1,...,N_t$, where
$\pi_{t 1},\ldots,\pi_{t N_t}$ is a feasible solution
of the respective dynamic program - program  \eqref{eq-6-1} for  $t=2,...,T-1$,  and program \eqref{eq-3a-1} for $t=T$.
We say that the dual problem \eqref{eq-2} has {\rm Relatively
Complete Recourse (RCR)} if    at every stage $t=2,...,T$,  for  any  generated $\pi_{t-1}$  by the forward process,    the respective dynamic program has a  feasible  solution at stage $t$ for every  realization of the random data.
\end{definition}

Without  RCR   it could  happen that   $V_t(\pi_{t-1},\xi_{t-1})=-\infty$ for a generated  $\pi_{t-1}$   and $\xi_{t-1}=\xi_{t-1 j}$. Unfortunately, it could happen that  the dual problem does not have the RCR property even if the primal problem has it. This could happen even in the two stage case.
One way to deal with     the  problem of absence of RCR  in numerical procedures is to use feasibility cuts, we will discuss this later. Another way  is the following penalty approach which will be used in Section \ref{sec:dualsddp}.
The infeasibility of problem \eqref{eq-3a-1} can happen because of its last constraint.  In order to deal with this, consider the following relaxation of problem \eqref{eq-3a-1}:
\begin{equation}\label{eq-r1}
\begin{array}{cll}
  \max\limits_{\pi_{T 1},...,\pi_{T N_T},\!\!,\,\zeta_T\ge 0} &    \sum\limits_{j=1}^{N_T}   p_{T j}  b_{T j}^\top \pi_{T j}-v_T^\top \zeta_T\\
    {\rm s.t.}&
  A_{T j}^\top \pi_{T j}  \le c_{T j},\;j=1,...,N_T,\\
    &A_{T -1}^\top \pi_{T-1}+
     \sum\limits_{j=1}^{N_T}   p_{T j}  B_{T j}^\top \pi_{T j} \le c_{T-1}+\zeta_T,
   \end{array}
\end{equation}
where $v_T$ is a  vector with positive components.
We have that problem \eqref{eq-r1} is always feasible
and hence its  optimal value $\tilde{V}_T(\pi_{T-1},\xi_{T-1})>-\infty$.
We also have that
\begin{equation}\label{valineq}
\tilde{V}_T(\pi_{T-1},\xi_{T-1})\ge V_T(\pi_{T-1},\xi_{T-1}),
\end{equation}
 with   the equality
 holding  if $\zeta_T=0$  in the optimal solution of \eqref{eq-r1}.  If
$V_T(\pi_{T-1}$,$\xi_{T-1})$ is finite, this equality holds if the components of vector $v_T$ are large enough.

Similarly, problems \eqref{eq-6-1} can be relaxed to
\begin{equation}\label{eq-r2}
\begin{array}{cll}
  \max\limits_{\pi_{t 1},\ldots,\pi_{t N_t}, \zeta_t \geq 0} &  \sum\limits_{j=1}^{N_{t}}  p_{t j}  \left[ b_{t j}^\top
  \pi_{t j}+ \tilde{V}_{t+1}(\pi_{t j},\xi_{t j})\right]-v_t^\top \zeta_t\\
    {\rm s.t.}&
  A_{t-1}^\top \pi_{t-1}+
  \sum\limits_{j=1}^{N_{t}} p_{tj} B_{t j}^\top  \pi_{t j}  \le c_{t-1}+\zeta_t,
   \end{array}
\end{equation}
with vector  $v_t$ having positive components.
In that way, the infeasibility problem is avoided and the  obtained  value gives an  upper bound for the optimal value of the dual problem.   Note that for sufficiently large vectors $v_t$ this upper bound coincides with the optimal value of the dual problem.


\subsection{Dynamics of Lagrange multipliers}
\label{sec-dyn}

Let us consider for the moment the two stage setting, i.e., $T=2$. The primal problem can be written as
\begin{equation}\label{lag-1}
    \min\limits_{x_1\ge 0}   c_1^\top x_1+ \bbe\left[ Q(x_1,\xi_2)\right]\;
    {\rm s.t.}\; A_1 x_1=b_1,
 \end{equation}
where $Q(x_1,\xi_2)$ is the optimal value of the second stage problem
\begin{equation}\label{lag-2}
    \min\limits_{x_2\ge 0}      c_2(\xi_2)^\top x_2\;
    {\rm s.t.}\;  B_2(\xi_2)  x_{1}+A_2(\xi_2)  x_2 =b_2(\xi_2).
\end{equation}
The Lagrangian of problem \eqref{lag-2} is
\[
L(x_1,x_2,\lambda,\xi_2)=c_2(\xi_2)^\top x_2+\lambda^\top (b_2(\xi_2)
- B_2(\xi_2)  x_{1}-A_2(\xi_2)  x_2).
\]
In the dual form,   $Q(x_1,\xi_{2 j})$ is given by the optimal value of the problem
 \begin{equation}\label{lag-3}
 \max\limits_{\lambda_j }  \; (b_{2 j}-B_{2 j} x_1)^\top \lambda_j   \;
  {\rm s.t.}
 \;  c_{2 j} -    A_{2 j}^\top \lambda_j   \ge 0.
\end{equation}
We have that if $x_1=\bar{x}_1$ is an optimal solution of the first stage problem, then optimal Lagrange multipliers $\pi_{2 j}$ are given by the optimal solution of problem \eqref{lag-3}.

This can be extended to the multistage setting of problem \eqref{eq-1} (recall that the stagewise independence condition is assumed). At the last stage $t=T$, given optimal solution $\bar{x}_{T-1}$, the following problem should be solved
\begin{equation}\label{lag-4}
    \min\limits_{x_T\ge 0}     c_T(\xi_T)^\top x_T \;\;
    {\rm s.t.}\;\;  B_T(\xi_T) \bar{x}_{T-1}+A_T(\xi_T) x_T=b_T(\xi_T).
\end{equation}
For a realization  $\xi_T=\xi_{T j}$,  the dual of problem \eqref{lag-4} is the problem
 \begin{equation}\label{lag-5}
 \max\limits_{\lambda_j } \;  (b_{T j}-B_{T j} \bar{x}_{T-1})^\top \lambda_j   \;\; {\rm s.t.}  \;\;  c_{T j} -  A_{T j}^\top \lambda_j   \ge 0.
\end{equation}
We then have that $\pi_{T j}$ are given by the optimal solution of problem \eqref{lag-5}.

At stage $t=T-1$, given optimal solution $\bar{x}_{T-2}$,  the following problem is supposed to  be solved (see  \eqref{primal})
\begin{equation}\label{lag-6}
\begin{array}{cll}
    \min\limits_{x_{T-1}\ge 0}  &   c_{T-1}(\xi_{T-1})^\top x_{T-1}+\Q_T(x_{T-1})\\
    {\rm s.t.} &  A_{T-1}(\xi_{T-1}) x_{T-1}=b_{T-1}(\xi_{T-1})-B_{T-1}(\xi_{T-1}) \bar{x}_{T-2}.
    \end{array}
\end{equation}
   We have that $\Q_T(\cdot)$  is a convex piecewise linear function. Therefore for every realization $\xi_{T-1}=\xi_{T-1 j}$  it is possible to represent \eqref{lag-6} as a linear program and hence to write its dual.  The optimal  Lagrange multipliers of that dual give the corresponding
Lagrange multipliers $\pi_{T-1 j}$.  And so on for other  stages going backward in time.  That is, we have the following.
\begin{remark}\label{remarkcomputedualmpsddp}  If $(\bar{x}_1,...,\bar{x}_T(\xi_{[T]}))$ is an optimal solution of the primal problem, then for $x_{t-1}=\bar{x}_{t-1}$ the Lagrange multiplier $\pi_{t j}$ is given by the respective  Lagrange multiplier of problem \eqref{primal}.
\end{remark}

We also refer to \cite{higlesenaor, rockdualsol} for the dynamics
of dual solutions to stochastic programs.

\setcounter{equation}{0}
\section{Sensitivity analysis}
\label{sec-sens}
In this section we discuss an application of the duality analysis to a study of sensitivity of the optimal value to small perturbations of the involved parameters.

\subsection{General case}
Suppose now that the data $c_t(\xi_t,\theta),  b_t(\xi_t,\theta), B_t(\xi_t,\theta)$, $A_t(\xi_t$, $\theta)$ of problem \eqref{eq-1} also depend  on parameter vector $\theta\in \bbr^k$. Denote by
$\vv(\theta)$  the optimal value of the parameterized problem   \eqref{eq-1}  considered as a function of $\theta$, and by $\cS(\theta)$ and $\cD(\theta)$ the sets of optimal solutions of the respective primal and dual problems. Recall that the sets $\cS(\theta)$ and $\cD(\theta)$ are nonempty provided the optimal value $\vv(\theta)$ is finite.
Let $L(x,\pi,\theta)$ be the corresponding Lagrangian (see \eqref{lagr}) considered  as a function of $\theta$.
Then
we have the following formula for the directional derivatives of the optimal value function (e.g.,  \cite[Proposition 4.27]{Bonnans2000}).

\begin{proposition}
\label{pr-valsens}
Suppose   that    the data functions are continuously differentiable functions of $\theta$, and   for a given $\theta=\bar{\theta}$ the optimal value $\vv(\bar{\theta})$ is finite and the sets $\cS(\bar{\theta})$ and $\cD(\bar{\theta})$ of optimal solutions   are bounded. Then
\begin{equation}\label{sen-1}
 \vv'(\bar{\theta},h)=\max_{\pi\in \cD(\bar{\theta})}\min_{x\in \cS(\bar{\theta})} h^\top \nabla_\theta L(x,\pi,\bar{\theta}).
\end{equation}
In particular if $\cS(\bar{\theta})=\{\bar{x}\}$ and $\cD(\bar{\theta})=\{\bar{\pi}\}$
are singletons, then $\vv(\cdot)$  is differentiable at $\bar{\theta}$ and
\begin{equation}\label{sen-2}
 \nabla \vv(\bar{\theta})= \nabla_\theta L(\bar{x},\bar{\pi},\bar{\theta}).
\end{equation}
\end{proposition}

Next, as an example, we consider the sensitivity analysis of an inventory model.

\subsection{Application to an inventory model}
\label{sec-invent}

Consider the inventory model
\begin{equation}\label{inv-1}
\begin{array}{cll}
 \min & \mathbb{E}\left[\displaystyle \sum_{t=1}^T
 a_t( y_t - x_{t-1} ) + g_t (\mathcal{D}_t - y_t)_{+} + h_t (y_t - \mathcal{D}_t)_{+} \right]\\
 {\rm s.t.}& x_{t} = y_t -\mathcal{D}_t, y_t \geq x_{t-1},t=1\ldots,T.
 \end{array}
\end{equation}
Here  $\mathcal{D}_1,...,\mathcal{D}_T$ is a (random) demand process,   $a_t,g_t,h_t$ are the ordering, back-order penalty and holding costs per unit, respectively,
$x_t$ is the inventory level and $y_t-x_{t-1}$ is the order quantity at time $t$,
 the initial inventory level $x_0$ is given. We refer to   \cite{Zipk} for a thorough discussion of that model.
Note that $\mathcal{D}_t$
is a random variable whereas
$d_t$ stands for a particular realization.
We assume that $g_t>a_t\ge 0$, $h_t> 0$, $t=1,...,T$.

In the classical setting the demand process is assumed to be stagewise independent, i.e., $\mathcal{D}_{t+1}$ is assumed to be independent of $\mathcal{D}_{[t]}=(\mathcal{D}_1,...,\mathcal{D}_t)$ for $t=1,...,T-1$.
In order to capture the autocorrelation structure of the demand process it is tempting  to model it as, say first order, autoregressive process
$
 \mathcal{D}_t=\mu+\phi \mathcal{D}_{t-1}  +\ep_t,
$
where errors  $\ep_t$ are assumed to be a sequence  i.i.d (independent identically distributed)  random variables. However this approach may result in some of the realizations of the demand process to be negative, which of course does not make sense. One way to deal with this is to make the transformation  $Y_t:=\log \mathcal{D}_t$ and to model $Y_t$  as an autoregressive process. A problem with this approach is  that it  leads to nonlinear equations for the original process $\mathcal{D}_t$,  which makes it difficult to use in the numerical algorithms discussed below.

We assume that
the demand  is modeled as the following multiplicative autoregressive process
\begin{equation}
\label{inv-2-1bis}
\mathcal{D}_t = \ep_{t}( \phi \mathcal{D}_{t-1} + \mu), \ t = 1,...,T,
\end{equation}
where $\phi\in (0,1)$, $\mu\ge 0$ are   parameters and $\mathcal{D}_0\ge 0$ is given. The errors $\ep_t$ are i.i.d  with  log-normal distributions having means and standard deviations given by $\bbe[\ep_t]=1$
and $\text{Var}(\ep_t)=\sigma^2>1$, respectively.
This guarantees that all realizations of the demand process are positive.   It is possible to view \eqref{inv-2-1bis} as a linearization of the log-transformed process
$\log \mathcal{D}_t$  (cf., \cite{sha2013}).
See Section \ref{statpr} for a discussion of statistical properties of the  process \eqref{inv-2-1bis}.

The process \eqref{inv-2-1bis} involves parameters $\phi$ and $\mu$ which are supposed to be estimated from the data. As such, these parameters are subject to estimation errors. This raises the question of sensitivity of the optimal value $\vv=\vv(\phi,\mu)$  of the corresponding problem \eqref{inv-1} viewed as a function of $\phi$ and $\mu$. To this end, we investigate the calculation of the derivatives
$\partial \vv(\phi,\mu)/\partial \phi$ and $\partial \vv(\phi,\mu)/\partial \mu$.  With these derivatives at hand, asymptotic distributions of the estimates of $\phi$ and $\mu$ can be translated into the asymptotics of the optimal value in a straightforward way by application of the Delta Theorem.  We refer to Section \ref{invnumer} for the corresponding numerical experiments.

\subsubsection{Properties of
 the multiplicative autoregressive process}
 \label{statpr}

Consider the multiplicative autoregressive process \eqref{inv-2-1bis}.
Note that under the specified conditions the demand process is not stationary. Indeed,
since the errors $\ep_t$ are i.i.d and $\bbe[\ep_t]=1$  we have  that
$
\bbe[\mathcal{D}_t ]=   \phi \bbe[\mathcal{D}_{t-1}] + \mu
$
and
\begin{equation}
\label{inv-3}
\begin{array}{lll}
\var(\mathcal{D}_t)&=&\bbe\left[\var\big ( \e_{t}( \phi \mathcal{D}_{t-1} + \mu)|\mathcal{D}_{t-1}\big )\right]+  \var  \left[ \bbe (\e_{t}( \phi \mathcal{D}_{t-1} + \mu)|\mathcal{D}_{t-1})\right] \\
&=&\bbe\left[ \sigma^2 ( \phi \mathcal{D}_{t-1} + \mu)^2\right]+\var( \phi \mathcal{D}_{t-1} + \mu)\\
&=&\sigma^2 \bbe\left[( \phi \mathcal{D}_{t-1} + \mu)^2\right]+\phi^2 \var(\mathcal{D}_{t-1} ).
\end{array}
\end{equation}
It follows that $\bbe[\mathcal{D}_t]$ converges to $\mu/(1-\phi)$ as $t\to \infty$.
 Suppose, for example,   that $\mu=0$. Then
 $
\mathcal{D}_t = \e_{t}  \phi \mathcal{D}_{t-1}=\mathcal{D}_0\phi^t\prod_{\tau=1}^t \e_\tau, \ t = 1,...,T,
$ $\bbe[\mathcal{D}_t]=\mathcal{D}_0 \phi^t\to 0$, and
 $
\var(\mathcal{D}_t )=\mathcal{D}_0^2\phi^{2t}[(1+\sigma^2)^t-1].
$
Therefore if $ \phi^2(1+\sigma^2)<1$, then
$\var(\mathcal{D}_t )\to 0$; and  if $ \phi^2(1+\sigma^2)> 1$, then
$\var(\mathcal{D}_t )\to \infty$ provided $\mathcal{D}_0>0$.

\section{Dual SDDP}\label{sec:dualsddp}

In this section, using the results of Section \ref{sec-dual}, we discuss  an adaptation of the cutting planes approach  for the approximation  of the value functions of the dual problem, similar to the standard SDDP method
and called Dual SDDP. The interested reader
can find the implementation
of  Primal SDDP and all variants of Dual SDDP described in this section at
{\url{https://github.com/vguigues/Dual_SDDP_Library_Matlab}} and {\url{https://github.com/vguigues/Primal_SDDP_Library_Matlab}}.

We will make the following assumption.
\begin{itemize}
\item[(A2)] Primal problem \eqref{eq-1}
satisfies the RCR assumption.
\end{itemize}

We first consider the case where only $b_t$
and $B_t$ are random in $\xi_t$.

\subsection{Dual SDDP for problems with uncertainty in $b_t$ and $B_t$}

In Dual SDDP, concave value functions $V_t,t=2,\ldots,T$, are approximated at the end of  iteration $k$ by polyhedral upper bounding functions
$V_t^k$ given by:
\begin{equation}\label{cutvtkdsddp}
V_t^k( \pi_{t-1}) = \displaystyle \min_{0 \leq i \leq k}
{\overline{\theta}}_t^i + \langle {\overline{\beta}}_t^i , \pi_{t-1} \rangle
\end{equation}
where
${\overline{\theta}}_t^i$, ${\overline{\beta}}_t^i$
are coefficients whose computation is detailed below. The algorithm uses valid upper bounds on the norm of dual optimal solutions:

\begin{lemma}\label{boundedpi}
Suppose that the optimal value of primal problem \eqref{eq-1} is finite  and that there is ${\hat x}>0$
feasible for primal problem
\eqref{eq-1}.
Then for every $t=1,\ldots,T$,
we can find
$\underline{\pi}_t, \overline{\pi}_t
\in \mathbb{R}^{m_t}$
such that dual problem
\eqref{eq-6-1}  is unchanged (i.e., has the
same optimal value) adding box
constraints
$\underline{\pi}_t \leq \pi_t \leq  \overline{\pi}_t$.
\end{lemma}

Recall that it is assumed that the number of scenarios is finite and hence problem \eqref{eq-1} can be viewed as a large linear program. The assumption of existence of feasible ${\hat x}>0$ means that   problem  \eqref{eq-1} possesses a feasible solution with all components being strictly positive.
If moreover  the equality constraints of problem \eqref{eq-1} are linearly independent, then this strict feasibility condition   implies that the set of optimal solutions of the dual problem (i.e., the set of
Lagrange multipliers) is bounded. On the other hand,  in the above lemma the linear independence condition is not assumed. A proof of Lemma \ref{boundedpi} and a way to obtain the corresponding
bounds $\underline{\pi}_t, \overline{\pi}_t$
can be found in the Appendix.

As mentioned earlier, a difficulty to solve the dual
problem with an SDDP type method is that
RCR may not be satisfied by the dual problem, even if RCR holds for the primal.
We propose two variants of Dual SDDP
to solve the Dual problem even if RCR does not hold for the dual:
Dual SDDP with penalizations and Dual SDDP with feasibility cuts.

\par {\textbf{Dual SDDP with penalizations.}} Dual SDDP with penalizations is based on the developments of
Section \ref{sec:rcr}. It introduces slack variables in the constraints which may become infeasible
for some past decisions in the subproblems solved in the forward passes of Dual SDDP.
Slack variables are penalized in the objective function with
sequences $(v_{t k})_k$  of positive penalizing coefficients.
Therefore, all subproblems solved in forward and backward passes of this variant of Dual SDDP, called
Dual SDDP with penalizations, are always feasible
and at iteration $k$, the method builds
polyhedral upper bounding function  $V_t^k$ for $V_t$ of form \eqref{cutvtkdsddp} (see Proposition \ref{proofubfeassddp}). Similarly to SDDP, trial points are generated in a forward pass and cuts for $V_t$
are computed in a backward pass. The detailed Dual SDDP method with penalizations
is as follows.\\

\par {\textbf{Initialization.}} For $t=2,\ldots,T,$ take for ${V}_t^0$
an affine upper bounding function for $V_t$
and $V_{T+1}^0 \equiv 0$. Set iteration counter $k$ to 1.\\

\par {\textbf{Step 1: forward pass of iteration $k$ (computation of dual trial points).}}
For the first stage of the forward pass, we compute an optimal solution
$\pi_1^k$ of
\begin{equation}\label{firststagepenalized}
V^{k-1} =
\begin{array}{l}
\displaystyle \max_{\pi_1} \; b_1^\top \pi_1 + V_{2}^{k-1}(\pi_1)\\
{\underline{\pi}}_1 \leq \pi_1 \leq {\overline{\pi}}_1.
\end{array}
\end{equation}
Recall that the optimal value of the first stage problem
does not change adding box constraints
${\underline{\pi}}_1 \leq \pi_1 \leq {\overline{\pi}}_1$  for appropriate values ${\underline{\pi}}_1$ and ${\overline{\pi}}_1$.
The introduction of these box constraints ensures that
the optimal value of \eqref{firststagepenalized} (which is an
{\em{approximate}} first stage problem
due to the approximation of
$V_2$ by $V_2^{k-1}$) is finite for all iterations.

\par For stage $t=2,\ldots,T-1$, given $\pi_{t-1}^k$, we compute an optimal solution of
\begin{equation}\label{forwarddualsddppen}
\begin{array}{cll}
  \max\limits_{\pi_{t 1},\ldots,\pi_{t N_t}, \zeta_t \geq 0} &  \sum\limits_{j=1}^{N_{t}}  p_{t j}  \left[ b_{t j}^\top
  \pi_{t j}+ {V}_{t+1}^{k-1}(\pi_{t j})\right]-v_{t k}^\top \zeta_t\\
    {\rm s.t.}&
  A_{t-1}^\top \pi_{t-1}^k +
  \sum\limits_{j=1}^{N_{t}} p_{tj} B_{t j}^\top  \pi_{t j}  \le c_{t-1}+\zeta_t,\\
  & {\underline{\pi}}_t \leq \pi_{t j} \leq {\overline{\pi}}_t.
   \end{array}
\end{equation}
\par An optimal solution of the problem above has $N_t$ components
$(\pi_{t 1}, \pi_{t 2},\ldots,\pi_{t N_t})$
for $\pi_t$. We generate a realization $\tilde \xi_t^k$ of
$\xi_t^k \sim \xi_t$ independently of previous realizations
$\tilde \xi_2^1,\ldots$, $\tilde \xi_{T-1}^1$,$\ldots$,
$\tilde \xi_2^k,\ldots,\tilde \xi_{t-1}^k$, and take $\pi_{t}^k = \pi_{t j_t(k)}$
where index $j_t(k)$ satisfies
$\tilde \xi_{t}^k = \xi_{t j_t(k)}$.\\

\par {\textbf{Step 2: backward pass of iteration $k$ (computation of new cuts)}}. We first compute a new cut for $V_T$. Let
$(\alpha, \delta, \overline{\Psi}, \underline{\Psi})$ be an optimal solution of\footnote{We suppressed the dependence of the optimal solution on $T$ and $k$ to alleviate notation.}
\begin{equation}\label{defVTk}
\begin{array}{l}
\displaystyle \min_{\alpha, \delta, \overline{\Psi}, {\underline{\Psi}}} \;
\delta^{\top} (c_{T-1} - A_{T-1}^{\top} \pi_{T-1}^k)  + c_T^{\top}  \sum_{j=1}^{N_T} \alpha_j
+ \sum_{j=1}^{N_T} {\overline{\Psi}}_j^\top {\overline \pi}_T
- \sum_{j=1}^{N_T} {\underline{\Psi}}_j^\top {\underline \pi}_T \\
A_T \alpha_j + p_{T j} B_{T j} \delta -
{\underline{\Psi}}_j +{\overline{\Psi}}_j = p_{T j} b_{T j},\;j=1,\ldots,N_T,\\
0 \leq \delta \leq v_{T k}, \alpha_j, {\underline{\Psi}}_j, {\overline{\Psi}}_j \geq 0\;j=1,\ldots,N_T.
\end{array}
\end{equation}
\par The new cut for $V_T$ has  coefficients
given by
$$
\overline \theta_T^k =
\delta^{\top} c_{T-1}+ c_T^{\top}  \sum_{j=1}^{N_T} \alpha_j
+ \sum_{j=1}^{N_T} {\overline{\Psi}}_j^\top {\overline \pi}_T
- \sum_{j=1}^{N_T} {\underline{\Psi}}_j^\top {\underline \pi}_T,\;
\overline \beta_T^k = -A_{T-1}\delta.
$$
\par For $t=T-1,\ldots,2$, compute an optimal solution
$(\delta,\nu,{\overline \Psi}, {\underline \Psi})$ of
\begin{equation}\label{defVtk}
\begin{array}{l}
\displaystyle \min_{\delta,\nu,{\overline \Psi}, {\underline \Psi}} \;
\delta^{\top}\Big[
c_{t-1} - A_{t-1}^{\top} \pi_{t-1}^k \Big]
+ \sum_{i=0}^{k} {\overline{\theta}}_{t+1}^i \sum_{j=1}^{N_t} \nu_i(j) +
\sum_{j=1}^{N_t} {\overline \Psi}_j^\top {\overline \pi}_t
-  \sum_{j=1}^{N_t} {\underline \Psi}_j^\top {\underline \pi}_t\\
\displaystyle
p_{t j} B_{t j} \delta - \sum_{i=0}^k \nu_i(j)  {\overline{\beta}}_{t+1}^i - {\underline \Psi}_j
+ {\overline \Psi}_j = p_{t j} b_{t j},\;j=1,\ldots,N_t,\\
\displaystyle \sum_{i=0}^k \nu_i( j ) = p_{t j}, {\underline \Psi}_j, {\overline \Psi}_j \geq 0,\;j=1,\ldots,N_t,\\
\nu_0,\ldots,\nu_k \geq 0,0 \leq \delta \leq v_{t k},
\end{array}
\end{equation}
\par and the cut coefficients
$$
\overline{\theta}_t^k = \delta^{\top} c_{t-1}
+ \sum_{i=0}^{k} {\overline{\theta}}_{t+1}^i \sum_{j=1}^{N_t} \nu_i(j) +
\sum_{j=1}^{N_t} {\overline \Psi}_j^\top {\overline \pi}_t
-  \sum_{j=1}^{N_t} {\underline \Psi}_j^\top {\underline \pi}_t,\;
{\overline{\beta}}_t^k = -A_{t-1} \delta.
$$
\par {\textbf{Step 3:}} Do $k \leftarrow k+1$ and go to Step 1.\\

\par
\begin{figure}
 \centering
 \begin{tabular}{c}
 \includegraphics[scale=0.5]{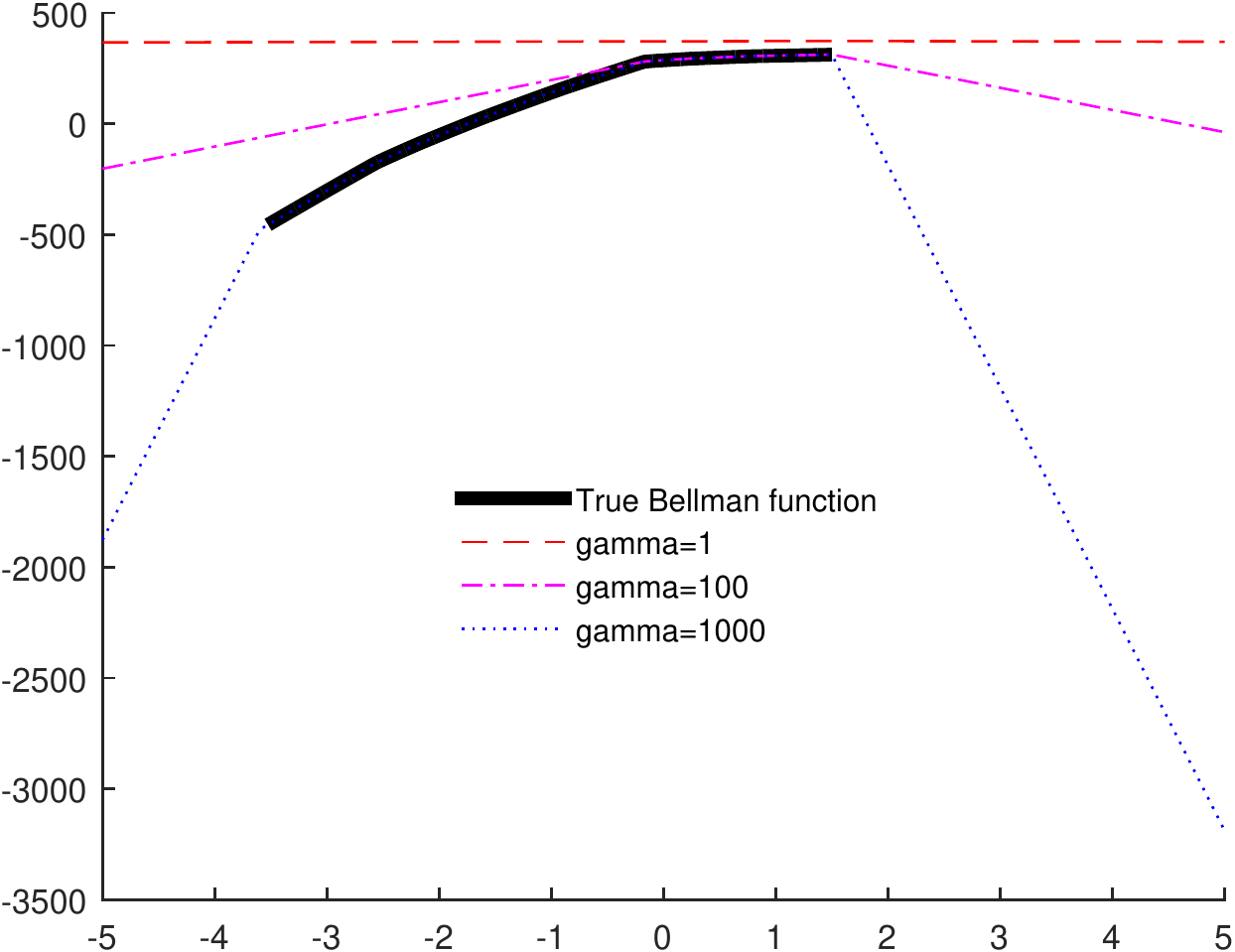}
 \end{tabular}
\caption{Graph of $V_2$ and of $V_2^{\gamma}$ for
$\gamma=1, 100, 1000$.}\label{fig1inventbisb}
\end{figure}

 The validity of the cuts computed in the backward pass of Dual SDDP with penalizations is shown in Proposition \ref{proofubfeassddp}.
\begin{proposition}\label{proofubfeassddp}
Consider Dual SDDP algorithm with penalizations.
Let Assumptions (A1) and (A2) hold.
Then for every $t=2,\ldots,T$, the sequence $V_t^k$ is a nonincreasing sequence
of upper bounding functions for $V_t$, i.e., for every $k \geq 1$ we have
$V_t \leq {V}_t^k \leq V_t^{k-1}$
 and therefore $(V^k)$ (recall that $V^{k-1}$ is the optimal value of \eqref{firststagepenalized}) is a nonincreasing deterministic
 sequence of upper bounds on the optimal value of \eqref{eq-1}.
\end{proposition}
To understand the effect of the sequence of penalizing parameters $(v_{t k})$
on Dual SDDP with penalizations, we define the following Dynamic Programming equations (see also
Lemma \ref{lem-conver} in the Appendix):
\begin{equation}\label{dualdpesimp1bb}
V_T^{\gamma}(\pi_{T-1}) =
\left\{
\begin{array}{cll}
  \max\limits_{\pi_{T 1},...,\pi_{T N_T},\!\!,\,\zeta_T\ge 0} &    \sum\limits_{j=1}^{N_T}   p_{T j}  b_{T j}^\top \pi_{T j}-\gamma {\textbf{e}}^\top \zeta_T\\
    {\rm s.t.}&
  A_{T j}^\top \pi_{T j}  \le c_{T j},\;j=1,...,N_T,\\
    &A_{T -1}^\top \pi_{T-1}+
     \sum\limits_{j=1}^{N_T}   p_{T j}  B_{T j}^\top \pi_{T j} \le c_{T-1}+\zeta_T,
   \end{array}
   \right.
   \end{equation}
for $t=2,\ldots,T-1$:
\begin{equation}\label{dualdpesimp2bb}
V_t^{\gamma}(\pi_{t-1}) =
\left\{
\begin{array}{cll}
  \max\limits_{\pi_{t 1},\ldots,\pi_{t N_t}, \zeta_t \geq 0} &  \sum\limits_{j=1}^{N_{t}}  p_{t j}  \left[ b_{t j}^\top
  \pi_{t j}+ V_{t+1}^{\gamma}(\pi_{t j})\right]-\gamma {\textbf{e}}^\top \zeta_t\\
    {\rm s.t.}&
  A_{t-1}^\top \pi_{t-1}+
  \sum\limits_{j=1}^{N_{t}} p_{tj} B_{t j}^\top  \pi_{t j}  \le c_{t-1}+\zeta_t,
   \end{array}
\right.
\end{equation}
and we define the first stage problem
 \begin{equation}\label{eq-7simpbb}
  \max_{\pi_1} \;\;      \pi_1^\top  b_1 +  V_{2}^{\gamma}(\pi_1),
\end{equation}
where {\textbf{e}} is a vector of ones and $\gamma$ is a positive real number.
As we will see below, $V_t^{\gamma}$ can be seen as an upper bounding concave approximation
of $V_t$ which gets ``closer'' to $V_t$ when $\gamma$ increases.
For inventory problem \eqref{inv-1}, it is easy to see that
functions $V_t$ in DP equations \eqref{eq-3a-1}, \eqref{eq-6-1}, \eqref{eq-7-1} and functions $V_t^{\gamma}$ in DP equations
\eqref{dualdpesimp1bb}, \eqref{dualdpesimp2bb}, \eqref{eq-7simpbb} (obtained
using in these equations data $c_t, b_t$ $A_t$, $B_t$, corresponding to the inventory problem) are only functions of one-dimensional state
variable $\pi_{t-1}$. Therefore, Dynamic Programming can be used to solve
these Dynamic Programming equations and obtain
good approximations of functions $V_t$ and $V_t^{\gamma}$. To obtain these approximations,
we need to obtain approximations of the domains of functions  $V_t$
and compute approximations of these functions on a set of points in that domain.
To observe the impact of penalizing term $\gamma$ on $V_t^{\gamma}$, we run Dynamic Programming both on DP equations \eqref{eq-3a-1}, \eqref{eq-6-1}, \eqref{eq-7-1}
and on DP equations \eqref{dualdpesimp1bb}, \eqref{dualdpesimp2bb}, \eqref{eq-7simpbb} for $\gamma=1$, $100$, and $1000$,
on  an instance of the inventory problem with $T=20$ and $N_t=20$.
The corresponding graphs of  $V_2$ (bold dark solid line) and of $V_2^{\gamma}$ for $\gamma=1, 10$, $1000$,
are represented in Figure \ref{fig1inventbisb}. We observe that all functions $V_2^{\gamma}$ are, as expected,
concave upper bounding functions for $V_2$ finite everywhere. We also see that on the domain of
$V_2$, $V_2^{\gamma}$ gets closer to $V_2$ when $\gamma$ increases and eventually coincides with $V_2$
on this domain when $\gamma$ is sufficiently large. Similar graphs were observed for remaining functions
$V_t, V_t^{\gamma}$, $t=3,\ldots,T$.
Therefore, convergence of Dual SDDP with penalizations requires
the coefficients $v_{t k}$ to become arbitrarily large. Proof of the following theorem is given in the Appendix.
\begin{theorem}\label{thconvdualsddppenal}
Consider optimization problem \eqref{eq-1} and Dual SDDP
with penalizations applied to the dual of this problem.
Let Assumptions (A1) and (A2) hold.
Assume that samples
$\xi_t^{\ell}$, $t=2,\ldots,T$, $\ell \geq 1$,
 in the forward passes are independent,
that $v_{t k+1}\geq v_{t k}$ for all $t,k$,
 and that
$\lim_{k \rightarrow +\infty}
v_{t k} = +\infty$ for all stage $t$. Then the sequence
$V^{k}$ is a deterministic sequence of upper bounds on the optimal value
of \eqref{eq-1} which converges almost surely to the optimal value of this problem.
\end{theorem}
\par {\textbf{Dual SDDP with feasibility cuts.}}
For dual problems not satisfying the RCR assumption, a subproblem for a given stage $t$
in the forward pass can be infeasible. In this situation, as was done in
Section 5 of \cite{guiguessiopt2016} for SDDP, we can build a feasibility cut
for stage $t-1$ and go back to the previous stage $t-1$ to resolve the problem
with that feasibility cut added, and so on until a sequence of feasible states
is obtained for all stages. In this context, no penalized slack
variables are used, neither in the forward nor in the backward pass.
Since the adaptations from \cite{guiguessiopt2016}
are simple, we skip the details of the derivations
of this SDDP method applied to the dual.
It will be tested in the numerical experiments of Section
\ref{numer}.

\if{

More precisely, we define for stage $t$ the set
$
\mathcal{S}_t:=\{ \pi_t: {\tilde C}_t^\top \pi_t   \leq {\tilde d}_t \}$
where matrices ${\tilde C}_t$ and vectors $\tilde d_t$ are updated along the iterations and such that
a feasible $\pi_t$ must belong to $\mathcal{S}_t$.
The detailed steps of the algorithm are given below.\\

\par {\textbf{Initialization.}} For $t=2,\ldots,T,$ take for $V_t^0$ an affine upper bounding function
for $V_t$. Set iteration counter $k$ to 1.\\
\par {\textbf{Step 1: forward pass (computation of dual trial points)}}.
\par We set $t=1$.
\par {\textbf{While }} $t\leq T$
\par $\hspace*{0.6cm}${\textbf{if }}$t=1$ we compute an optimal solution $\pi_1^k$ of
$$
V^{k-1} =
\begin{array}{l}
\displaystyle \max_{\pi_1} \; b_1^\top \pi_1 + V_{2}^{k-1}(\pi_1)\\
{\underline{\pi}}_1 \leq \pi_1 \leq {\overline{\pi}}_1,
\end{array}
$$
\par $\hspace*{1.2cm}$where $\underline{\pi}_1$, $\overline{\pi}_1$ are as in Lemma \ref{boundedpi} and do $t \leftarrow t+1$.
\par $\hspace*{0.6cm}${\textbf{else }}

\par $\hspace*{1.2cm}${\textbf{if }}$t=T$
\par $\hspace*{1.4cm}$//Defining function $\mathcal{V}_{T}$ by
\begin{equation}\label{forwardfeas0}
\mathcal{V}_{T}(\pi_{T-1}):=
\begin{array}{l}
  \min \limits_{\pi_{T 1},\ldots,\pi_{T N_T}, \zeta_T} \;
{\textbf{e}}^\top \zeta_T \\
A_T^{\top} \pi_{T j} \leq c_T,\;j=1,\ldots,N_T,\\
\zeta_T \geq 0,
A_{T-1}^\top \pi_{T-1} +
\sum_{j=1}^{N_T} p_{T j} B_{T j}^{\top} \pi_{T j}
\leq c_{T-1}+\zeta_T,
  \end{array}
\end{equation}
\par $\hspace*{1.4cm}$//where as before $\textbf{e}$ is a vector of ones then if $\mathcal{V}_{T}(\pi_{T-1}^k)>0$ we have
\par $\hspace*{1.4cm}$//that $\pi_{T-1}^k$ is not feasible
for stage $T-1$.
\par $\hspace*{1.4cm}$//For any $\pi_{T-1}$, problem \eqref{forwardfeas0} is feasible with finite optimal value
\par $\hspace*{1.4cm}$//$\mathcal{V}_{T}(\pi_{T-1})$
which can be expressed as the optimal value of the
\par $\hspace*{1.4cm}$//corresponding dual problem, i.e.,
\begin{equation}\label{forwardfeas1}
\mathcal{V}_{T}(\pi_{T-1}):=
\begin{array}{l}
  \max \limits_{\delta, \alpha} \; -\sum\limits_{j=1}^{N_{T}}
 c_T^{\top} \alpha_j + \delta^\top (-c_{T-1}+A_{T-1}^{\top} \pi_{T-1})\\
\alpha_j \geq 0, p_{T j} B_{T j} \delta + A_T \alpha_j  = 0,\;j=1,\ldots,N_T,\\
0 \leq \delta \leq e.
  \end{array}
\end{equation}
\par $\hspace*{1.8cm}$Let $\delta,\alpha$ be an optimal
solution of \eqref{forwardfeas1}
written for $\pi_{T-1}=\pi_{T-1}^k$
\par $\hspace*{1.8cm}$with optimal value $\mathcal{V}_{T}(\pi_{T-1}^k)$.
\par $\hspace*{1.8cm}${\textbf{if }}$\mathcal{V}_{T}(\pi_{T-1}^k)>0$ then we
add a feasibility cut adding
to ${\tilde C}_{T-1}^{\top}$ the
\par $\hspace*{1.8cm}$row $(A_{T-1} \delta)^{\top}$ and to $\tilde d_{T-1}$
the component $\delta^\top c_{T-1} +
c_T^\top \sum\limits_{j=1}^{N_{t}} \alpha_j$.
We
\par $\hspace*{1.8cm}$then go back to the previous stage setting $t \leftarrow t-1$. The corresponding
\par $\hspace*{1.8cm}$feasibility cut
is valid since for $\pi_{T-1}$ feasible for stage
$T-1$ we have
$$\mathcal{V}_{T}(\pi_{T-1})=0 \geq -\sum\limits_{j=1}^{N_{T}}
 c_T^{\top} \alpha_j + \delta^\top (-c_{T-1}+A_{T-1}^{\top} \pi_{T-1}).$$

\par $\hspace*{1.8cm}${\textbf{else }}$t \leftarrow t+1$
\par $\hspace*{1.8cm}${\textbf{end if}}
\par $\hspace*{1.2cm}${\textbf{else }}
\par $\hspace*{1.8cm}$//Let $\mathcal{V}_t$ be the function given by
\begin{equation}\label{forwardfeas2p}
\mathcal{V}_t(\pi_{t-1}):=\begin{array}{l}
  \min \limits_{\zeta_t \geq 0,\pi_{t 1},\ldots,\pi_{t N_t}} \; {\textbf{e}}^{\top} \zeta_t\\
A_{t-1}^{\top} \pi_{t-1}
+ \sum_{j=1}^{N_t} p_{t j} B_{t j}^{\top} \pi_{t j}
\leq c_{t-1}+ \zeta_t,\\
{\tilde C}_t^{\top} \pi_{t j} \leq  \tilde d_t, j=1,\ldots,N_t.
  \end{array}
\end{equation}
\par $\hspace*{1.8cm}$//Linear program \eqref{forwardfeas2p} is feasible
and has a finite optimal value for
\par $\hspace*{1.8cm}$//every $\pi_{t-1}$. Therefore
by linear programming duality we can express
\par $\hspace*{1.8cm}$//$\mathcal{V}_t(\pi_{t-1})$
as the optimal value of the corresponding dual problem \par $\hspace*{1.8cm}$//given by
\begin{equation}\label{forwardfeas2}
\mathcal{V}_t(\pi_{t-1}):=\begin{array}{l}
  \max \limits_{\delta, \lambda} \; -\sum\limits_{j=1}^{N_{t}}
 {\tilde d}_t^{\top} \lambda_j +
 \delta^\top (-c_{t-1} + A_{t-1}^{\top} \pi_{t-1})\\
p_{t j} B_{t j} \delta + {\tilde C}_t {\lambda}_j = 0,j=1,\ldots,N_t,\\
0 \leq \delta \leq {\textbf{e}}, \lambda_j \geq 0,\;j=1,\ldots,N_t.
  \end{array}
\end{equation}
\par $\hspace*{1.8cm}$We compute an optimal solution
$\delta, \lambda$ of \eqref{forwardfeas2} written
with
\par $\hspace*{1.8cm}$$\pi_{t-1}=\pi_{t-1}^k$ with optimal value
$\mathcal{V}_t(\pi_{t-1}^k)$.
\par $\hspace*{1.8cm}${\textbf{if }}$\mathcal{V}_{t}(\pi_{t-1}^k)>0$ then we add
add to ${\tilde C}_{t-1}^{\top}$ the row $(A_{t-1} \delta)^{\top}$, we add
\par $\hspace*{1.8cm}$to $\tilde d_{t-1}$
the component $\delta^\top c_{t-1}+\sum\limits_{j=1}^{N_{t}} {\tilde d}_t^{\top} \lambda_j$, and do $t \leftarrow t-1$. The
\par $\hspace*{1.8cm}$corresponding
feasibility cut
is valid since for $\pi_{t-1}$ feasible for stage
\par $\hspace*{1.8cm}$$t-1$ we have
$$\mathcal{V}_{t}(\pi_{t-1})=0 \geq -\sum\limits_{j=1}^{N_{t}}
 {\tilde d}_t^{\top} \lambda_j +
 \delta^\top (-c_{t-1} + A_{t-1}^{\top} \pi_{t-1}).$$
\par $\hspace*{1.8cm}${\textbf{else }}we compute an optimal solution $\pi_{t}^k$ of
\begin{equation}\label{forwardfeas3}
\begin{array}{cll}
  \max\limits_{\pi_{t 1},\ldots,\pi_{t N_t}} &  \sum\limits_{j=1}^{N_{t}}  p_{t j}  \left[ b_{t j}^\top
  \pi_{t j}+ {V}_{t+1}^{k-1}(\pi_{t j})\right]\\
    {\rm s.t.}&
  A_{t-1}^\top \pi_{t-1}^k +
  \sum\limits_{j=1}^{N_{t}} p_{tj} B_{t j}^\top  \pi_{t j}  \le c_{t-1},\\
  & {\underline{\pi}}_t \leq \pi_{t j} \leq {\overline{\pi}}_t, {\tilde C}_t^{\top} \pi_{t j} \leq \tilde d_t,\;j=1,\ldots,N_t.
   \end{array}
\end{equation}
\par $\hspace*{1.8cm}$An optimal solution of the problem above has $N_t$ components
$(\pi_{t 1}, \pi_{t 2}$,
\par $\hspace*{1.8cm}$ $\ldots,\pi_{t N_t})$
for $\pi_t$. We generate a realization $\tilde \xi_t^k$ of
$\xi_t^k \sim \xi_t$
\par $\hspace*{1.8cm}$independently of previous realizations
$\tilde \xi_2^1,\ldots$, $\tilde \xi_{T-1}^1$,$\ldots$,
$\tilde \xi_2^k,\ldots,\tilde \xi_{t-1}^k$,
\par $\hspace*{1.8cm}$and take $\pi_{t}^k = \pi_{t j_t(k)}$ where index $j_t(k)$ satisfies
$\tilde \xi_{t}^k = \xi_{t j_t(k)}$.
\par $\hspace*{1.8cm}$We then go to the next
stage: $t \leftarrow t+1$.
\par $\hspace*{1.8cm}${\textbf{end if}}
\par $\hspace*{1.2cm}${\textbf{end if}}
\par $\hspace*{0.6cm}${\textbf{end if}}
\par {\textbf{end while}}\\
\par {\textbf{Step 2: backward pass}}. For
$t=2,\ldots,T$, we compute a new cut for $V_t$ to build $V_t^k$ exactly as in Dual SDDP
with penalizations taking penalizations
$v_{t k}=0$ (which is possible since trial points
are feasible).\\
\par {\textbf{Step 3:}} Do $k \leftarrow k+1$ and go to Step 1.\\

In Dual SDDP with feasibility cuts
all subproblems solved in the backward passes are feasible and $V^k$ is a nonincreasing deterministic
 sequence of upper bounds on the optimal value of \eqref{eq-1}.
Moreover, following the proof of
Theorem 5.1 in \cite{guiguessiopt2016}, we obtain
that if the samples in the forward passes
are independent then
the sequence
$V^{k}$ converges to the optimal value of
primal problem \eqref{eq-1}.

}\fi

\subsection{Dual SDDP for problems with uncertainty in all parameters}

We have seen in Section \ref{wrdpequ} how to write DP equations
on the dual problem of a MSLP when all data $(A_t,B_t,c_t,b_t)$ in $(\xi_t)$
is random. In this situation, cost-to-go functions $V_t$ are functions
$V_t(\pi_{t-1},\xi_{t-1})$ of both past decision $\pi_{t-1}$
and past value $\xi_{t-1}$ of process $(\xi_t)$.
Also recall that functions
$V_t(\cdot,\xi_{t-1})$ are concave for all $\xi_{t-1}$.
Therefore, Dual SDDP with penalizations from the previous section
must be modified as follows. For each stage $t=2,\ldots,T,$ instead of computing
just one approximation of a single function (function $V_t$), we now need
to compute approximations of $N_t$ functions, namely
concave cost-to-go functions $V_t(\cdot,\xi_{t-1 j})$, $j=1,\ldots,N_t$.
The approximation $V_{t j}^k$ computed for $V_t(\cdot,\xi_{t-1 j})$  at iteration $k$ is a polyhedral function
$V_{t j}^k$ given by:
$$
V_{t j}^k( \pi_{t-1}) = \displaystyle \min_{0 \leq i \leq k}  \;{\overline{\theta}}_{t j}^i + \langle {\overline{\beta}}_{t j}^i , \pi_{t-1} \rangle.
$$
Therefore more computational effort is needed.
However, the adaptations of the method can be easily written.
More specifically, at iteration $k$, in the forward pass, dual trial points
are obtained replacing $V_{t}(\cdot,\xi_{t-1 j})$ by $V_{t j}^{k-1}$
and in the backward pass a cut is computed at stage $t$ for $V_{t}(\cdot,\xi_{t-1 j_k})$ with
$j_k$ satisfying $\xi_{t-1 j_k}=\tilde \xi_{t-1}^k$ where
$\tilde \xi_{t-1}^k$ is the sampled value of $\xi_{t-1}$ at iteration $k$.

\subsection{Dual SDDP for problems with interstage dependent cost coefficients}
\label{dualsddp}

We consider problems of form
\eqref{eq-1} where costs  $c_t$ affinely depend on their  past
while  $b_t$  are stagewise  independent.
Specifically, similar to derivations of Section \ref{sec-invent},
suppose that $c_t$ follow  a multiplicative vector autoregressive process of form
\begin{equation}\label{inv-2}
\begin{array}{ll}
 c_t = \varepsilon_t \circ \left(  \sum_{j=1}^p \Phi_{t j} c_{t-j} + \mu_t  \right),
 \end{array}
\end{equation}
with  $(x \circ y)_i=x_i y_i$ denoting the componentwise product,  and where matrices $\Phi_{t j}$ and vectors  $\mu_t \geq 0$ as well as $c_1,\ldots,c_{2-p} \geq 0$ are given.

We assume that the  process
$(b_t,\varepsilon_t)$ is stagewise
independent and that the support
of $b_t, \varepsilon_t$ is the finite
set
$$\{
(b_{t 1},\varepsilon_{t 1}),
\ldots,
(b_{t N_t},\varepsilon_{t N_t})
\},$$
with
$\varepsilon_{t i}>0$
and  $p_{t i}=
\mathbb{P}\{
(b_t, \varepsilon_t)=
(b_{t i}, \varepsilon_{t i})\}$,
$i=1,\ldots,N_t$.
For some values of $\Phi_{t j}$ (for instance for matrices with
nonnegative entries), this guarantees that all realizations of the price process $\{c_t\}$ are positive.
The developments which follow
can be easily extended to other linear models for $\{c_t\}$, for instance SARIMA or PAR
models, see \cite{vgcoap2014} for the definition of state vectors of minimal
size for generalized linear models.

Using the notation
$c_{t_1:t_2}=(c_{t_1},c_{t_1+1},\ldots,c_{t_2-1},c_{t_2})$ for
$t_1 \leq t_2$ integer, for the corresponding primal problem (of the form \eqref{eq-1}),  we can write the following Dynamic Programming equations:
define $\mathcal{Q}_{T+1} \equiv 0$ and for $t=2,\ldots,T$,
\begin{equation} \label{defqtb}
\mathcal{Q}_t ( x_{t-1} , c_{t-p:t-1} ) =
\mathbb{E}_{b_t, \varepsilon_t}\Big[
Q_t ( x_{t-1} , c_{t-p:t-1}, b_t, \varepsilon_t ) \Big]
\end{equation}
where $Q_t (x_{t-1} , c_{t-p:t-1}, b_t, \varepsilon_t )$ is given by
\begin{equation}\label{deffraktb}
\begin{array}{ll}
\displaystyle \min_{x_t \geq 0}  \Big[\varepsilon_t \circ \Big(\displaystyle \sum_{j=1}^p \Phi_{t j} c_{t-j} + \mu_t \Big) \Big]^{\top}
x_t + \mathcal{Q}_{t+1}\Big(x_t,
c_{t+1-p:t-1},\varepsilon_t \circ \Big(\displaystyle \sum_{j=1}^p \Phi_{t j} c_{t-j} + \mu_t \Big)\Big ) \\
 A_t x_t + B_t x_{t-1} = b_t,
\end{array}
\end{equation}
while the first stage problem is
$$
\begin{array}{l}
\displaystyle \min_{x_1 \geq 0} \; c_1^\top x_1 + \mathcal{Q}_2( x_1,c_{2-p:1}) \\
A_1 x_1 = b_1.
\end{array}
$$

Standard  SDDP does not apply directly to solve Dynamic Programming equations
\eqref{defqtb}-\eqref{deffraktb} because functions $\mathcal{Q}_t$ given by \eqref{defqtb}-\eqref{deffraktb} are not convex.
Nevertheless,   we can use the Markov Chain discretization variant of SDDP to solve
Dynamic Programming equations \eqref{defqtb}-\eqref{deffraktb}.
On the other hand, as pointed above,  it is possible to apply SDDP for the dual problem with the added state variables.
Along the lines of Section \ref{wrdpequ}
we can write Dynamic Programming equations for the dual,
now with function $V_t$ depending on
$\pi_{t-1},c_{t-1},\ldots,c_{t-p}$.

\if{
The following Dynamic Programming equations
for the dual of \eqref{eq-1} with $(c_t)$ of form \eqref{inv-2} can be be written.
For the last stage $T$, we have to solve the problem:
\begin{equation}\label{dualdpesimp1b}
\begin{array}{l}
  \max\limits_{\pi_{T 1},\ldots,\pi_{T N_T}} \; \sum\limits_{j=1}^{N_T} p_{T j} \pi_{T j}^{\top} b_{T j}\\
    A_{T}^{\top} \pi_{T j} \leq
    \varepsilon_{T j} \circ \left(\mu_{T} +  \sum_{\ell=1}^p  \Phi_{T \ell} c_{T-\ell} \right),\;j=1,...,N_T,\\
      \sum_{j=1}^{N_T} p_{T j} B_{T}^{\top} \pi_{T j} \leq c_{T-1}-A_{T-1}^{\top} \pi_{T-1},
   \end{array}
\end{equation}
with optimal value $V_T(\pi_{T-1} , c_{T-p:T-1})$.

Next for stage $t=2,\ldots,T-1$, given $V_{t+1}$, we need
to solve the problem
\begin{equation}\label{dualdpesimp2b}
\begin{array}{cl}
 \max \limits_{\pi_{t 1},\ldots,\pi_{t N_t}} & \sum \limits_{j=1}^{N_{t}}  p_{t j}
\left( \pi_{t j}^\top b_{t j} + V_{t+1}\Big(\pi_{t j},  c_{t+1-p:t-1},\varepsilon_{t j} \circ (
\mu_{t} +  \sum_{\ell=1}^p  \Phi_{t \ell} c_{t-\ell} )\Big) \right)\\
{\rm s.t.} &  \sum_{j=1}^{N_{t}} p_{t j}
  B_{t}^{\top} \pi_{t j}  \leq c_{t-1} - A_{t-1}^{\top} \pi_{t-1},
   \end{array}
   \end{equation}
while the first stage problem is
 \begin{equation}\label{eq-7simpb}
  \max_{\pi_1} \;\;      \pi_1^\top  b_1  +  V_{2}(\pi_1 , c_{2-p:1}).
\end{equation}
}\fi
These functions are concave
and  therefore we can apply
Dual SDDP with penalizations to these DP equations
to build polyhedral approximations of these functions $V_t$ of form
\begin{equation}
V_t^k(\pi_{t-1},c_{t-1},\ldots,c_{t-p})
= \min_{0 \leq i \leq k} \theta_{t}^i +  \langle \beta_{t 0}^i , \pi_{t-1} \rangle + \sum_{j=1}^p  \langle \beta_{t j}^i , c_{t-j} \rangle
\end{equation}
at iteration $k$.

\setcounter{equation}{0}
\section{Numerical experiments}
\label{numer}

In this section, we report numerical results obtained applying Primal SDDP and  variants of Dual SDDP
to the  inventory problem and to the  Brazilian interconnected power system problem.
All methods were implemented in Matlab and run on an Intel Core i7, 1.8GHz, processor
with 12,0 Go of RAM. Optimization problems were solved using Mosek \cite{mosek}.
 \if{
When the data process  of problem \eqref{eq-1} is  interstage dependent, in some cases the problem can be reformulated  in a stagewise independent form by adding  appropriate state variables. In case the data process can   represented as an autoregressive process and only the  right hand sides $b_t$  of \eqref{eq-1} are random, this preserves the linearity of the obtained problem (see  section \ref{sec-invent}). On the other hand if the costs $c_t$ are given by an
autoregressive process, then such approach of adding state variables destroys the linearity and convexity of the primal problem. It could be noted that  in the dual problem \eqref{eq-2} the costs are in the right hand side  of the constraints,  and   hence such approach can be applied to the dual problem.
}\fi

\subsection{Dual SDDP for the inventory problem}\label{sec:inventorysim}

\if{
namely as
$$
-a_1 x_0 +
\left\{
\begin{array}{l}
\displaystyle \min \mathbb{E}_{\xi}\left[  \sum_{t=1}^T c_t^{\top} x_t + d_t^{\top} u_t  \right] \\
A_1 x_1 + C_1 u_1 = b_1, D_1 x_1 + F_1 u_1 \leq  f_1,\\
A_t x_t + B_t x_{t-1} + C_t u_t = b_t, \mbox{a.s.},\;t=2,\ldots,T,\\
D_t x_t + E_t x_{t-1} + F_t u_t \leq  f_t,\mbox{a.s.},\;t=2,\ldots,T,
\end{array}
\right.
$$
where the underlying stochastic process $\mathcal{D}_t$ is the demand process, control $u_t=[y_t;w_t;v_t]$,
$$
\begin{array}{lllll}
c_t=  \left\{
\begin{array}{ll}
-a_{t+1} & \mbox{if }t<T,\\
0 & \mbox{for }t=T,
\end{array}
\right.
&
A_1= [-1],
&
C_1=[1 \;  0  \; 0],
&
D_1=[0;0;0;0;0].
&
F_1=\left[
\begin{array}{ccc}
0 & -1 & 0 \\
-1 & -1 & 0 \\
0 & 0 & -1 \\
1 & 0 & -1 \\
-1 & 0 & 0
\end{array}
\right],
\end{array}
$$
$$
\begin{array}{llll}
d_t=\left[
\begin{array}{c}
a_t \\
g_t \\
h_t
\end{array}
\right],
b_t=[\mathcal{D}_t], &
f_1= \left[
\begin{array}{c}
0\\
-\xi_1\\
0\\
\xi_1 \\
-x_0
\end{array}
\right],
 &   f_t= \left[
\begin{array}{c}
0\\
-\mathcal{D}_t\\
0\\
\mathcal{D}_t \\
0
\end{array}
\right], t \geq 2,
\end{array}
$$
and for $t=2,\ldots,T$:
$$
\begin{array}{llllll}
A_t=[-1],&B_t=[0],&C_t=[1 \; 0 \; 0],&D_t=[0;0;0;0;0]&E_t=[0;0;0;0;1],&F_t=F_1.
\end{array}
$$
}\fi
We consider the inventory problem \eqref{inv-1}
with parameters $a_t=1.5+\cos(\frac{\pi t}{6})$, $p_{t i}=\frac{1}{N}$ where $N$ is the number of realizations for each stage, $\xi_{t j}=(5+0.5 t)(1.5+0.1 z_{t j})$ where
$(z_{t 1},\ldots,z_{t N})$ is a sample from the standard Gaussian
distribution, $x_0=10$, $g_t=2.8$, and
$h_t=0.2$.

\par {\textbf{Illustrating the correctness of DP equations \eqref{eq-3a-1}, \eqref{eq-6-1}, \eqref{eq-7-1} and checking the convergence of the variants of Dual SDDP.}}
\begin{figure}
 \centering
 \begin{tabular}{cc}
  \includegraphics[scale=0.48]{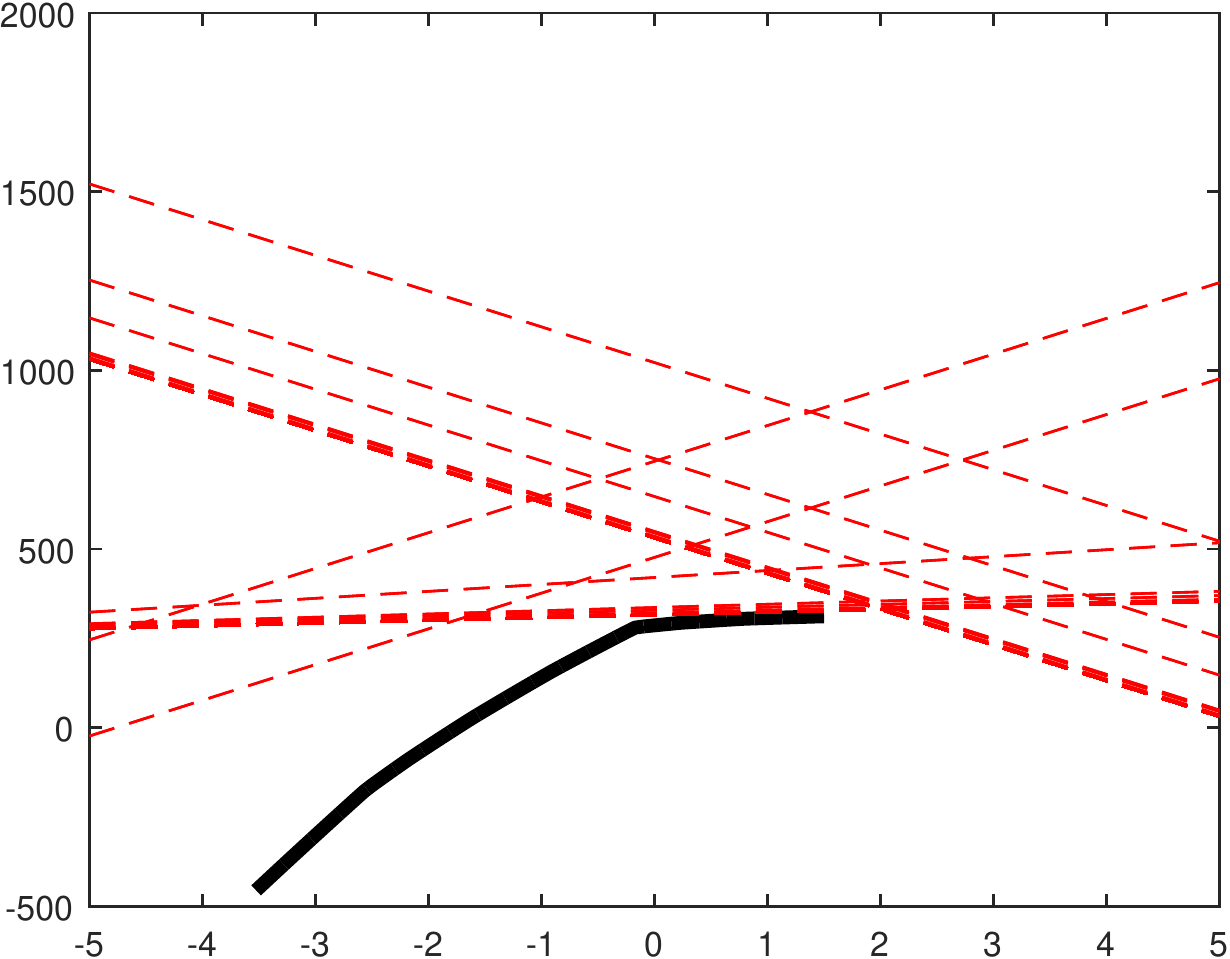} &
\includegraphics[scale=0.48]{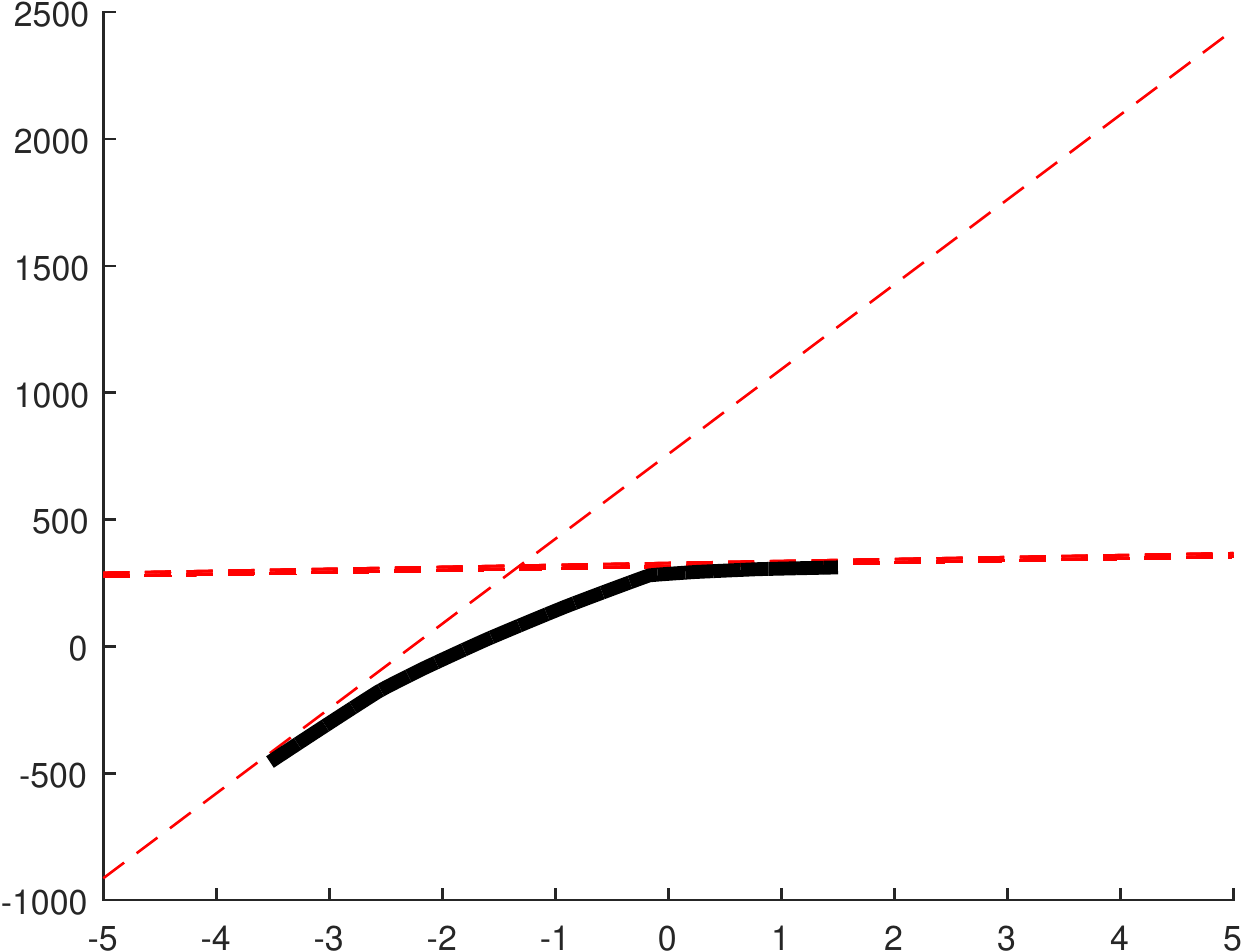}
 \end{tabular}
 \caption{Graph of $V_2$ (bold black solid line) and cuts computed for $V_2$
by Dual SDDP with penalizations $v_{t k}=100$ (left panel) and
Dual SDDP with feasibility cuts (right panel).}\label{fig1invent}
\end{figure}
We solve this inventory problem using Dynamic Programming applied both
to DP equations \eqref{eq-3a-1}, \eqref{eq-6-1}, \eqref{eq-7-1} and to DP equations \eqref{dualdpesimp1bb}-\eqref{dualdpesimp2bb}
for $\gamma=1, 10, 1000$.
In this latter case, we obtain approximations of functions $V_t^{\gamma}$.
We also run Primal SDDP, Dual SDDP with feasibility cuts,
and Dual SDDP with penalties $v_{t k}=1$, $10$, $1000$,
on the same instance, knowing that Dual SDDP variants were run
for 100 iterations (the upper bounds computed by these methods
stabilize in less than 10 iterations) and Primal SDDP was stopped
when
the gap is $<0.1$ where the gap is defined as $\frac{Ub-Lb}{Ub}$ where $Ub$ and $Lb$ correspond to upper and lower bounds
 computed
by Primal SDDP along iterations.
The lower bound $Lb$ is the optimal value of the first stage problem and the
upper  bound $Ub$ is the upper end of a 97.5\%-one-sided confidence interval on the optimal value obtained using the sample
of total costs computed by all previous forward passes. With this stopping criterion and the considered instance
of the inventory problem, Primal SDDP was run for 232 iterations.

In Figure \ref{fig1invent}, we report the graph of $V_2$ and the cuts computed for
$V_2$ by Dual SDDP with feasibility cuts (right panel) and
Dual SDDP with penalties $v_{t k}=100$ (left panel).
All cuts are, as expected, upper bounding affine functions for $V_2$
on its domain. However, it is interesting to notice that for Dual SDDP with feasibility
cuts, few different cuts are computed and these cuts are tangent or very close to $V_2$
at the trial points.
On the contrary, Dual SDDP with penalties
may compute many cuts dominated by others on the domain of $V_2$.
Therefore, cut selection techniques, for instance along the lines of
\cite{guiguesejor17} \cite{guiguesbandarra18} using Limited Memory Level 1 cut selection, could be useful for Dual SDDP.

We report in Table \ref{table0}
the approximate optimal values  and the time needed to compute them
with Primal SDDP, Dual SDDP, and Dynamic Programming
applied to respectively \eqref{eq-3a-1}, \eqref{eq-6-1}, \eqref{eq-7-1}
and \eqref{dualdpesimp1bb}, \eqref{dualdpesimp2bb}, \eqref{eq-7simpbb} with $\gamma=1, 100, 1000$.
The approximate optimal values reported are the last upper bound computed for variants of Dual SDDP
and the last lower bound computed for Primal SDDP. All approximate optimal values
are very close (showing that all variants were correctly implemented) and Dynamic Programming is much slower than the other sampling-based algorithms.
For Dual SDDP with penalization, if penalties are too small the upper bound can be $+\infty$ while if penalties
are sufficiently large the algorithm converges to an optimal policy.

\begin{figure}
 \centering
 \begin{tabular}{cc}
\includegraphics[scale=0.48]{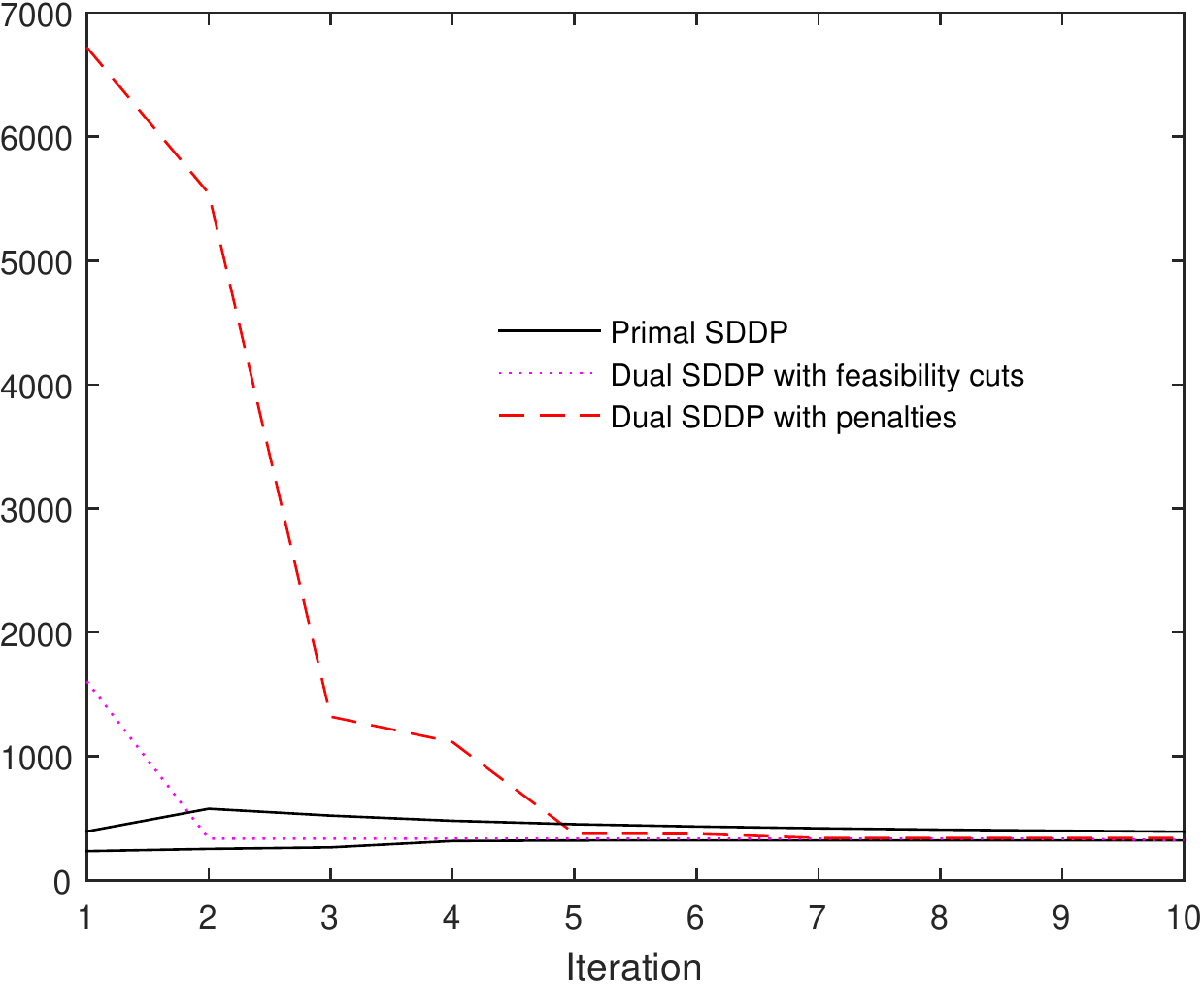}
 &
\includegraphics[scale=0.48]{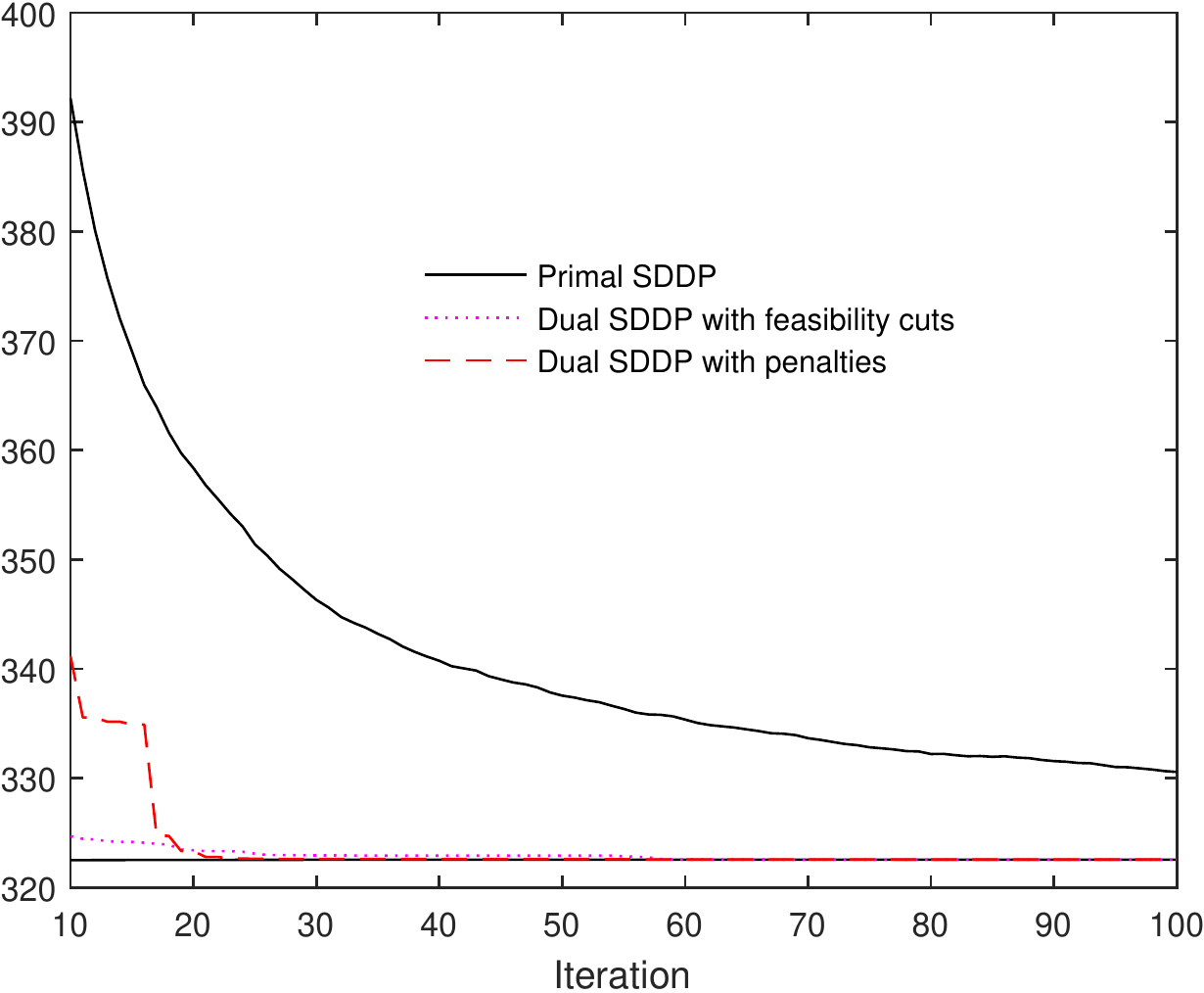}\\
 \end{tabular}
\caption{Left: upper and lower bounds computed by Primal SDDP
and upper bounds computed by Dual SDDP with feasibility cuts
and Dual SDDP with penalties $v_{t k}=1000$
for the first 10 iterations.
Right: same outputs for iterations 10,$\ldots,100$.}\label{fig1bounds}
\end{figure}

\begin{table}
\centering
\begin{tabular}{|c|c|c|}
\hline
{\tt{Method}} & {\tt{Optimal value}} &  {\tt{CPU time (s.)}}    \\
\hline
{\tt{DP}} on \eqref{eq-3a-1}, \eqref{eq-6-1}, \eqref{eq-7-1}  &   321.6& 685\\
\hline
{\tt{DP}} on  \eqref{dualdpesimp1bb}, \eqref{dualdpesimp2bb}, \eqref{eq-7simpbb}, $\gamma=1$  & $+\infty$ & 2 860 \\
\hline
{\tt{DP}} on  \eqref{dualdpesimp1bb}, \eqref{dualdpesimp2bb}, \eqref{eq-7simpbb}, $\gamma=100$  &  322.2 & 3 808\\
\hline
{\tt{DP}} on \eqref{dualdpesimp1bb}, \eqref{dualdpesimp2bb}, \eqref{eq-7simpbb}, $\gamma=1000$  & 321.8 & 3 376 \\
\hline
{\tt{Primal SDDP}}   &  322.5  &     105     \\
\hline
{\tt{Dual SDDP with penalties}}, $v_{t k}=1$  & 2 131.4 & 9.4\\
\hline
{\tt{Dual SDDP with penalties}}, $v_{t k}=100$  & 322.5& 11.3\\
\hline
{\tt{Dual SDDP with penalties}}, $v_{t k}=1000$  &  322.5 & 11.9\\
\hline
{\tt{Dual SDDP with feasibility cuts}}  &  322.5 &10.6 \\
\hline
\end{tabular}
\caption{Optimal value and CPU time needed (in seconds) to compute them on an instance of the inventory problem
with $T=N_t=20$ by Dynamic Programming (DP), Primal SDDP, and variants of Dual SDDP.}\label{table0}
\end{table}
Finally, we report for this instance in Figure \ref{fig1bounds}
the evolution of the lower bound $Lb$ and upper bound $Ub$ computed by Primal SDDP and the upper bounds computed by
Dual SDDP with penalties $v_{t k}=1000$ and Dual SDDP with feasibility cuts.
With  Dual SDDP, the upper bound is naturally large at the first iteration but decreases much quicker than the upper
bound $Ub$ computed by Primal SDDP, especially for Dual SDDP with feasibility cuts, with all upper bounds
converging to the optimal value of the problem.\\

\par {\textbf{Tests on a larger instance.}} We now run Primal and Dual SDDP on a larger instance
with $T=100$ and $N_t=100$ for 600 iterations.
The evolution of the upper bounds computed along the iterations
of Dual SDDP (both with feasibility cuts and with penalizations $v_{t k}=1000$)
and of the upper and lower bounds computed by Primal SDDP
are reported in Table \ref{tablebounditer}
for
iterations $2, 3, 5, 10, 50, 100, 200, 300, 400, 500,$ and $600$.
We see that for the first iterations, the upper bound decreases more quickly with the variants
of Dual SDDP, the most important decrease being obtained for Dual SDDP with feasibility cuts.
However, on this instance, the convergence of Dual SDDP with feasibility cuts is slower, i.e.,
a solution of high accuracy is obtained quicker using Dual SDDP with  penalizations.
More precisely, we fix confidence levels $\varepsilon=0.2, 0.15, 0.1, 0.05, 0.01$, and
for each confidence level, we compute the time needed, running Primal and Dual SDDP in parallel, to obtain a solution with relative accuracy $\varepsilon$ stopping the algorithm when
the upper bound {\tt{Ub\_D}} computed by a variant of Dual SDDP and the lower bound Lb, computed by Primal SDDP,
satisfies ({\tt{Ub\_D}}-Lb)/{\tt{Ub\_D}}$<\varepsilon$. The results are reported in Table \ref{table1}.
In this table, we also report the time needed to obtain a solution of relative accuracy
$\varepsilon$ using only the information provided by Primal SDDP, stopping the algorithm when
(Ub-Lb)/Ub$<\varepsilon$.

We observe that if $\varepsilon$ is not too small, the smallest CPU time
is obtained combining Primal SDDP with Dual SDDP
with feasibility cuts while when $\varepsilon$ is small (0.05 and 0.01) the smallest CPU
time is obtained combining Primal SDDP with Dual SDDP with penalizations. For $\varepsilon=0.05$ and $0.01$,
600 iterations are even not enough to get a solution of relative accuracy $\varepsilon$ using
Primal SDDP or combining Primal SDDP and Dual SDDP with feasibility cuts.

\if{
\begin{figure}
{\color{blue} this figure takes the whole page, can we remove it?}
 \centering
 \begin{tabular}{cc}
\includegraphics[scale=0.48]{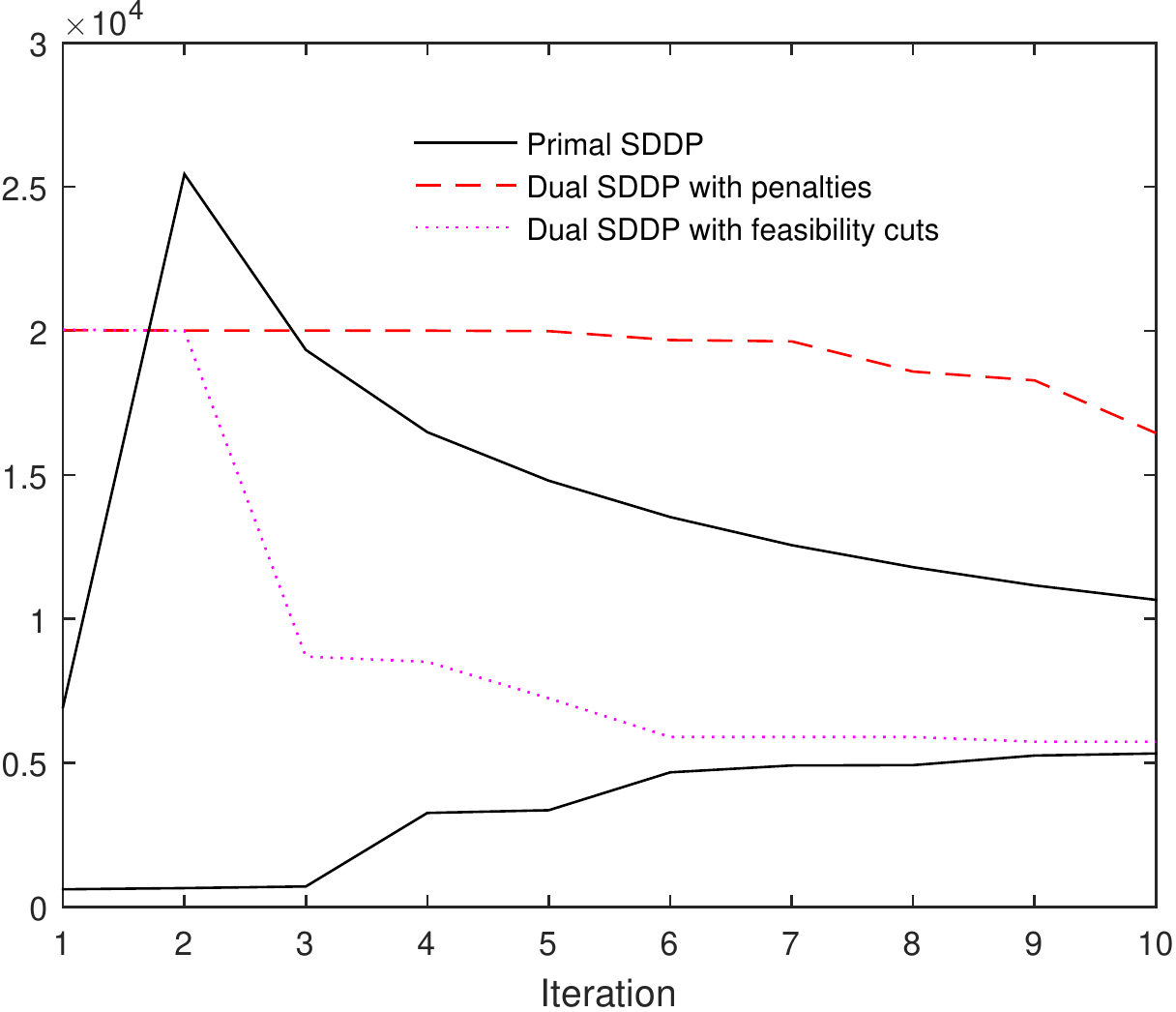}
 &
\includegraphics[scale=0.48]{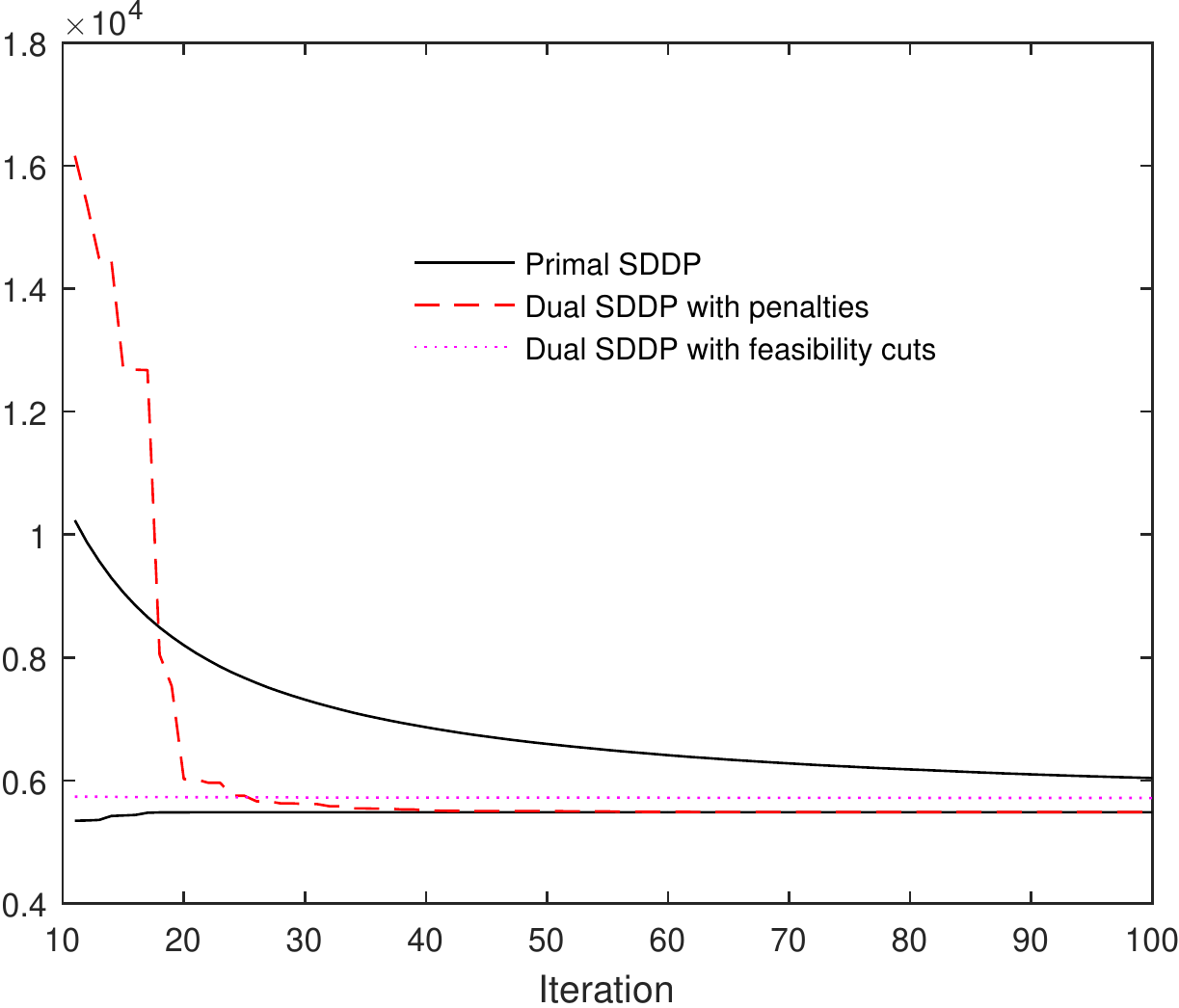}\\
 \end{tabular}
\caption{Left: upper and lower bounds computed by Primal SDDP
and upper bounds computed by Dual SDDP with feasibility cuts
and Dual SDDP with penalties $v_{t k}=1000$
for the first 10 iterations for an instance of the inventory problem with $T=N_t=100$.
Right: same outputs for iterations 10,$\ldots,100$.}\label{fig2bounds}
\end{figure}
}\fi

\begin{table}[H]
\centering
\begin{tabular}{|c|c|c|c|c|}
\hline
Iteration & \begin{tabular}{c}{\tt{Primal}}\\{\tt{SDDP}}\\Lb\end{tabular} & \begin{tabular}{c}{\tt{Primal}}\\
{\tt{SDDP}}\\Ub\end{tabular}  &  \begin{tabular}{c}{\tt{Dual SDDP with}}\\{\tt{feasibility}}\\{\tt{cuts}}\end{tabular} & \begin{tabular}{c}{\tt{Dual SDDP}}\\{\tt{with}}\\  {\tt{penalties}}\end{tabular}\\
\hline
2  & 656.4 &25 443  & 20 002 &20 015 \\
\hline
3  & 713.1  &19 340  &8 693.1  & 20 012\\
\hline
5  & 3361.8  &  14 800  &7 246.8&19 993\\
\hline
10 & 5330.1   &  10 662 &  5 736.6 & 16 452 \\
\hline
50     &    5483.1    &6 594.5& 5721.8 & 5500.9\\
\hline
100     &  5483.5     & 6 039.2&  5715.1 & 5484.8\\
\hline
200     &  5483.6     &5 762.4 & 5710.0  & 5484.2\\
\hline
300     &   5483.7   & 5 671.0 &  5704.6 & 5484.0\\
\hline
400     &  5483.7    &5 625.3  & 5702.7 & 5483.9\\
\hline
500     &  5483.7    &5 597.9  &  5702.5 & 5483.8\\
\hline
600     &   5483.7   & 5 579.9 &  5702.2 &5483.8 \\
\hline
\end{tabular}
\caption{For an instance of the inventory problem with $T=N_t=100$,
lower bound Lb and upper bound Ub computed by Primal SDDP and upper bounds computed
by Dual SDDP with
feasibility cuts and Dual SDDP with
penalties $v_{t k}=1000$ along iterations.}\label{tablebounditer}
\end{table}

\begin{figure}
 \centering
\includegraphics[scale=0.6]{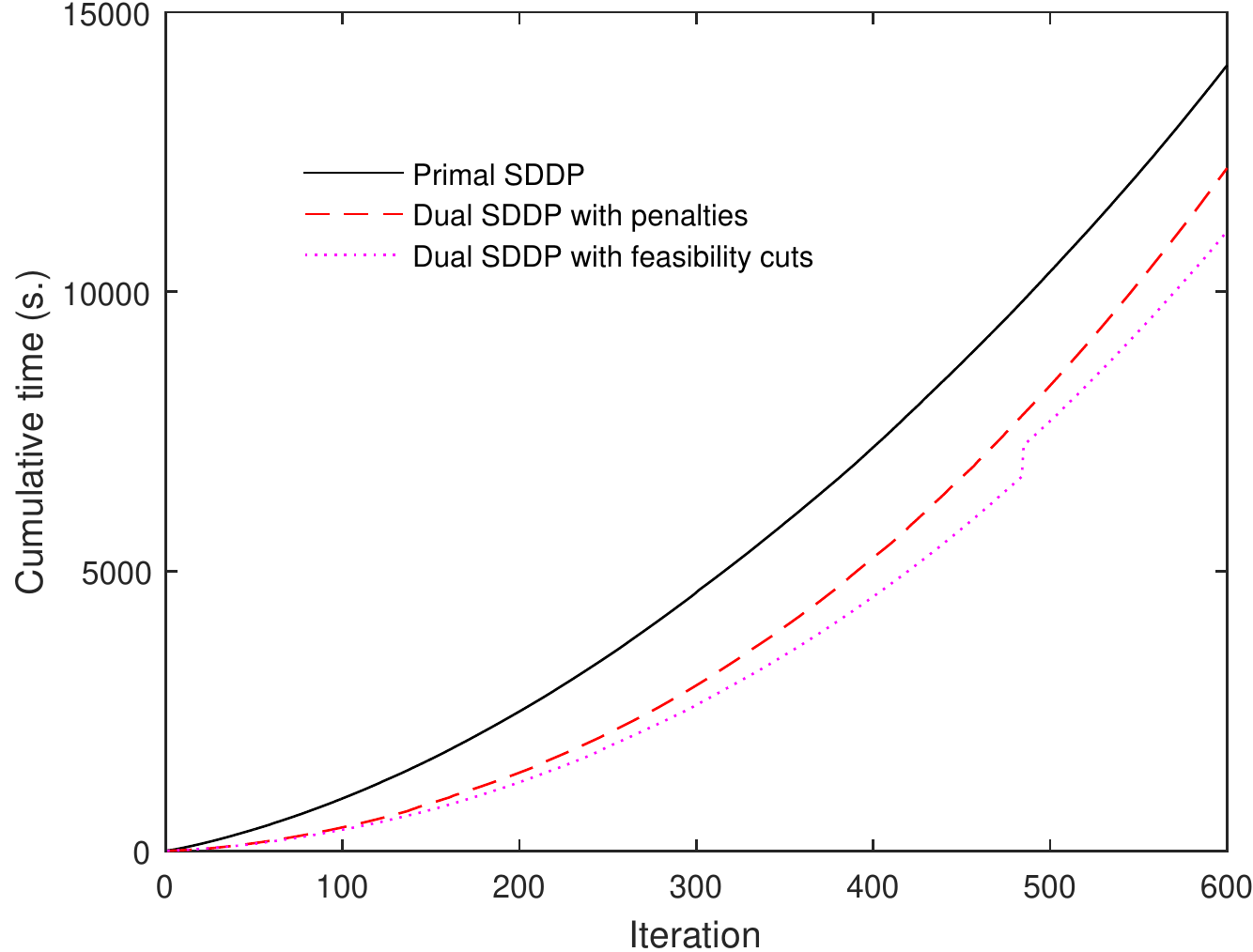}
\caption{Cumulative CPU time along iterations of Primal SDDP, Dual SDDP with feasibility cuts, and
Dual SDDP with penalizations $v_{t k}=1000$.}\label{figtimeinvent}
\end{figure}

In Figure \ref{figtimeinvent}, we report the cumulative CPU  time along iterations of
 all methods. We see that each iteration requires a similar computational bulk and the CPU time
 increases exponentially with the number of iterations.
 \begin{table}
\centering
\begin{tabular}{|c|c|c|c|}
\hline
$\varepsilon$ & {\tt{Primal SDDP}} & \begin{tabular}{c}{\tt{Dual SDDP with}}\\{\tt{feasibility cuts}}\end{tabular} & \begin{tabular}{c}{\tt{Dual SDDP with}}\\{\tt{penalties $v_{t k}=1000$}}\end{tabular}\\
\hline
0.2  & 300.2 & 29.5 & 35.8 \\
\hline
0.15  & 459.8  & 35.8 & 41.2 \\
\hline
0.1  & 825.6  & 48.3   &48.3\\
\hline
0.05 & 2366.2   &96.1&  61.5  \\
\hline
0.01     &    -    & -  &103.2 \\
\hline
\end{tabular}\\
\caption{Time needed (in seconds) to obtain a solution of relative accuracy $\varepsilon$
with Primal SDDP, Dual SDDP with feasibility cuts, and Dual SDDP with penalties $v_{t k}=1000$
for an instance of the inventory problem with $T=N_t=100$.}\label{table1}
\end{table}

\if{
Finally, for Dual SDDP with penalizations, we report in Figure \ref{fig3error} the maximal and mean values of
$\zeta_t^{k}$ (recall that for the inventory problem $\zeta_t$ in \eqref{forwarddualsddppen}
are real-valued)  along iterations for each stage $t$, where $\zeta_t^{k}$ is an optimal
value
of $\zeta_t$ in \eqref{forwarddualsddppen} for iteration $k$. The corresponding values
are positive, meaning that indeed RCR does not hold for the dual of the inventory problem.
}\fi

\if{
\begin{figure}
 \centering
 \begin{tabular}{cc}
\includegraphics[scale=0.4]{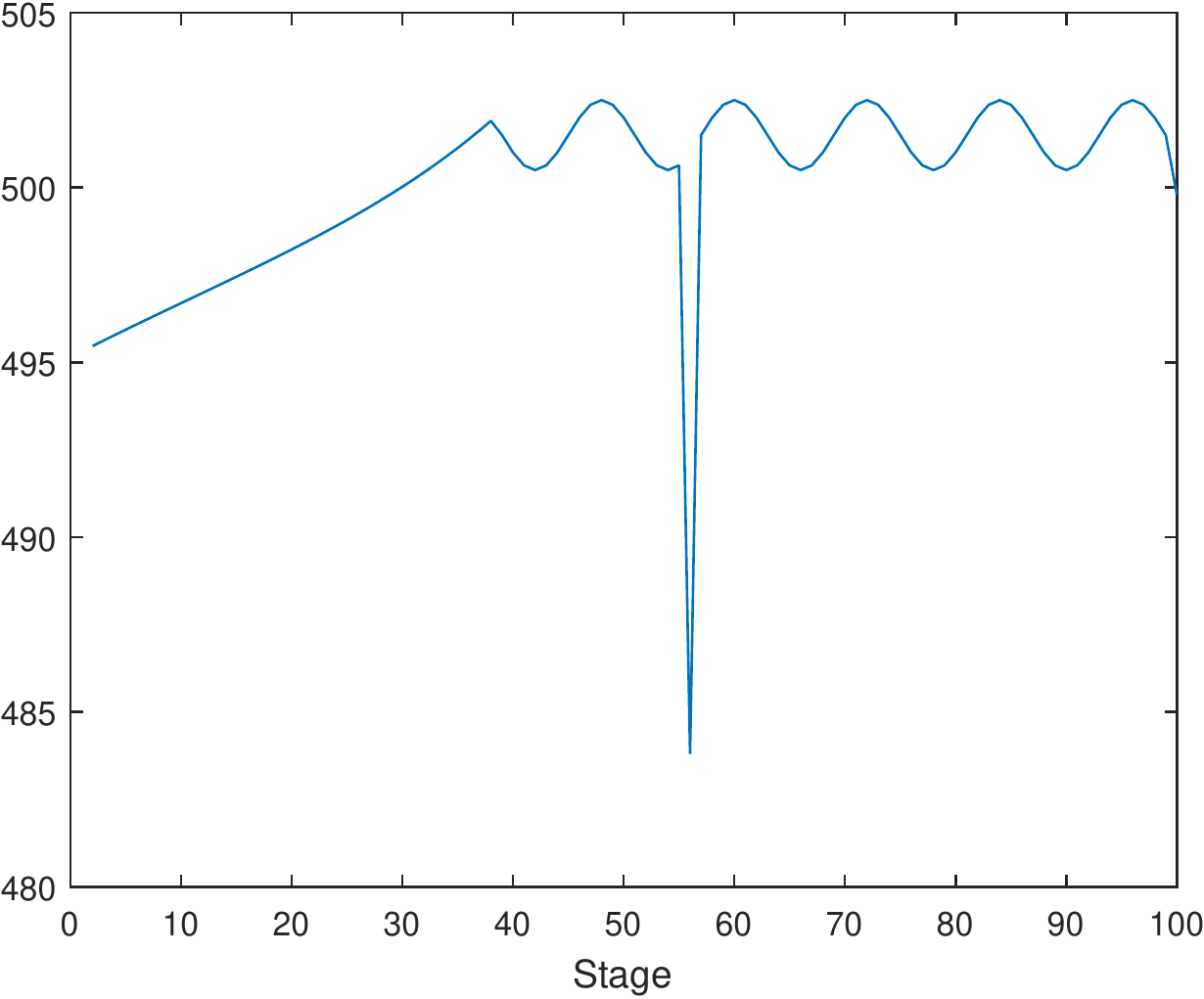}
 &
\includegraphics[scale=0.4]{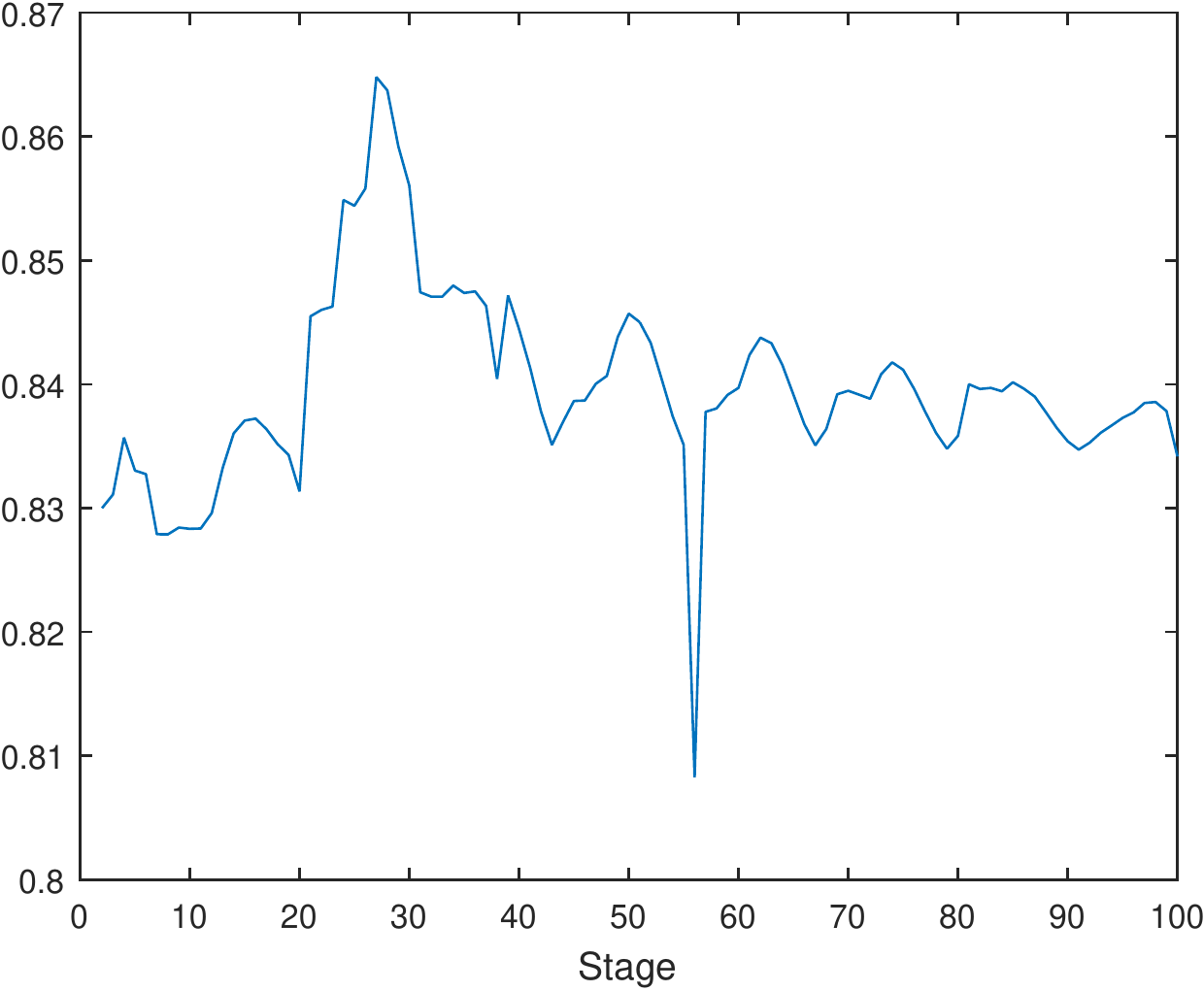}\\
 \end{tabular}
\caption{Left: maximal constraint violation as a function of the stage. Right: Mean constraint violation as a function of
the stage. }\label{fig3error}
\end{figure}
}\fi

\subsection{Sensitivity analysis for the inventory problem}
\label{invnumer}

Consider the inventory problem of Section \ref{sec:inventorysim} with $(\mathcal{D}_t)$ as in
\eqref{inv-2-1bis} and $T=10$ stages.
For this problem, the derivatives
from Proposition \ref{pr-valsens}
are given by
\begin{eqnarray}
\label{inv-12}
&\partial \vartheta(\phi,\mu)/\partial \phi=
\partial L(\bar{x},\bar{y},\bar{\pi})/\partial \phi=
\bbe\left[\sum_{t=1}^T
\bar{ \pi}_t  \epsilon_{t}  \mathcal{D}_{t-1}
 \right] ,\\
 \label{inv-13}
& \partial \vartheta(\phi,\mu)/\partial \mu=
 \partial L(\bar{x},\bar{y},\bar{\pi})/\partial \mu=
\bbe\left[\sum_{t=1}^T
\bar{ \pi}_t  \epsilon_{t}
 \right] ,
\end{eqnarray}
where $(\bar{x},\bar{y})$ is an optimal solution of the primal  problem and $\bar{\pi}$ are the corresponding Lagrange multipliers.
Our goal is to compute these derivatives
 solving
the primal and dual problems by
respectively Primal and Dual SDDP.

We consider $4$ instances with $(\phi,\mu)=(0.01,0.1)$, $(0.01,3.0)$, $(0.001,0.1)$, and $(0.001,3.0)$.
The remaining parameters of these instances are those from the previous section.
We discretize both the primal and dual problem into $N_t = 100$ samples
for each stage $t = 2,\ldots,10$. We take the relative error
$\varepsilon = 0.01$ for the stopping criterion and use
$10\,000$ Monte Carlo simulations to estimate the expectations in \eqref{inv-12}, \eqref{inv-13}.
For Primal SDDP, the upper bound Ub and lower bound Lb at termination are given
in Table $\ref{tab-sen-num-primal}$ for the four instances.
\begin{table}[htbp]
	\centering
	\small
	\begin{tabular}{|c|c|c|c|c|}
		\hline
	Bound		 & {\texttt{Instance 1}} &
			{\texttt{Instance 2}} & {\texttt{Instance 3}} & {\texttt{Instance 4}} \\
			\hline
	   Ub &17.9176    &  478.687 & 15.3940 & 404.242 \\
	   \hline
	    Lb	&17.9163    &475.017  & 15.3927 & 402.913 \\
	    \hline
	\end{tabular}
	\caption{Upper and lower bounds at the last iteration of Primal SDDP.}
	\label{tab-sen-num-primal}
\end{table}

The optimal mean values of Lagrangian multipliers for the demand constraints computed, for a given stage $t \geq 2$, averaging
over the $10\,000$ values obtained simulating $10\,000$ forward passes after termination, are given in Table \ref{tab-sen-num-dual}. In this table, \texttt{LM} stands for the multipliers obtained using Primal SDDP as explained in Remark \ref{remarkcomputedualmpsddp}
whereas \texttt{Dual} stands for the multipliers obtained using Dual SDDP with penalties.
The fact that the multipliers obtained are close for both methods illustrates the validity of the two alternatives we
discussed in Sections \ref{sec-sens}-\ref{sec:dualsddp} to compute derivatives of the value function of a MSP.
\begin{table}[htbp]
	\centering
	\small
	\begin{tabular}{|c|c|c|c|c|c|c|c|c|}
\hline
		Stage & \multicolumn{2}{|c|}{\texttt{Instance 1}} & \multicolumn{2}{|c|}{\texttt{Instance 2}} & \multicolumn{2}{|c|}{\texttt{Instance 3}} & \multicolumn{2}{|c|}{\texttt{Instance 4}} \\
		\hline
		{} &\texttt{LM} & \texttt{Dual} &\texttt{LM} & \texttt{Dual} &\texttt{LM} & \texttt{Dual} &\texttt{LM} & \texttt{Dual}\\
    \hline
        2	 &0.2465
     & 0.2373  &1.6701	&1.66959 & 0.0444	&0.0328 &1.666	&1.666\\
     \hline
        3    &0.3218      &  0.31095 &1.4098	&1.4120 &0.1421	&0.1340 &1.406	&1.409\\
         \hline
        4     & 0.3268   &	0.3221
  &0.9862	&0.9861 &0.19439	&0.18974 &0.984 &0.984\\
   \hline
        5     &0.3086   &	0.3058 & 0.6330	&0.6329 &0.2145	&0.2128  &0.6327	&0.6327\\
         \hline
        6     &0.3408   &	0.3412 &0.49998	& 0.499897   &0.2708	 &0.2717  &0.4999	&0.4998\\
         \hline
        7     &0.5026 	&0.5051   &0.63397	&0.63397   &0.4378	&0.4418 &0.6339	&0.6339\\
         \hline
       	8    &0.7047    & 0.7049
 &0.8348	& 0.8340 &0.6404 &	0.6413 &0.8349	&0.8334\\
  \hline
        9     &0.8985  &	0.9032 &1.0322	&1.0343    &0.83501	&0.8401 &1.0315	&1.0343\\
         \hline
       	10   &1.1022	& 1.1037 &1.2302 &1.2365  & 1.03926	&1.04091 &1.23	&1.23\\
       	 \hline
	\end{tabular}
	\caption{Comparison between optimal Lagrange multipliers from Primal SDDP and Dual SDDP with penalties.}
	\label{tab-sen-num-dual}
\end{table}

With optimal dual solutions $\{\bar \pi_t\}$ and the realizations of $\{D_t\}$ and $\{\epsilon_t\}$ at hand, we are able to compute the sensitivity of the optimal value with respect to $\phi$ and $\mu$, using (\ref{inv-12}) and (\ref{inv-13}),
with expectations estimated for $10\,000$ Monte Carlo simulations. We benchmark our method against the finite-difference  method. Specifically, for value function
$\vv$, the finite-difference method approximates the derivative with respect to $u_0$ by
$
	v'(u_0 ) \approx \frac{v(u_0 + \delta) - v(u_0 -\delta)}{2\delta}
	\label{sen-num:inv-3}$ for some small $\delta$.
	
	The sensitivity of the optimal value of the inventory problem with respect to $(\phi,\mu)$ is displayed in
	Table \ref{tab-sen-num-sensitivity}. In this table, S-$\phi$ and S-$\mu$ denote the derivatives with respect to
	$\phi$ and $\mu$ computed by our method, and fd-$\phi$, fd-$\mu$ denote the derivatives computed by the finite-difference method. In order to measure the difference between the two methods, we also compute S-gap-$\phi$ and S-gap-$\mu$, where
	S-gap-$\phi := \frac{|\text{fd-}\phi - \text{S-}\phi |}{|\text{fd-}\phi|} \times 100\%$ and S-gap-$\mu := \frac{|\text{fd-}\mu - \text{S-}\mu |}{|\text{fd-}\mu|} \times 100\%$.

 \begin{table}[htp]
 	\centering
 	\small
 	\begin{tabular}{|c|c|c|c|c|c|c|}
\hline
 		{\texttt{Instance}} & {\texttt{fd-}$\phi$} & {\texttt{S-}$\phi$} & {\texttt{S-gap-}$\phi$($\%$)}  & {\texttt{fd-}$\mu$} & {\texttt{S-}$\mu$} & {\texttt{S-gap-}$\mu$($\%$)}\\
 		\hline
 		1  & 403.604 & 401.094 & 0.622  & 164.578 & 164.158 & 0.255 \\
 		\hline
 		2   & 10\,716.111 & 10\,671.262 & 0.419 & 185.346 & 184.847 & 0.270 \\
 		\hline
 		3  & 269.514 & 269.443 & 0.026 & 134.646 & 134.463 & 0.136 \\
 		\hline
 		4   & 7\,780.570 & 7\,770.274 & 0.132 & 158.017 & 158.001 & 0.0101\\
 		\hline
 	\end{tabular}
 	\caption{Sensitivity of the optimal value with respect to $\phi$ and $\mu$ by the two methods.}
 	\label{tab-sen-num-sensitivity}
 \end{table}
We observe that the derivatives obtained by both methods are close to each other, especially when $\phi$ and $\mu$ are small. This is because small $\phi$ and $\mu$ gives rise to less variability in the demand. Note also that  the
finite-difference method is more time consuming since it requires computing the optimal value twice.
Instead, our method only needs to solve the model once. Moreover, computing the Lagrange multipliers does not significantly consume CPU time, as they are generated as a by-product of Primal SDDP. Alternatively, as discussed above, one can compute the optimal multipliers using Dual SDDP with penalties. Another drawback of the finite-difference method lies in its numerical instability. Indeed, the method is more accurate when $\delta$ is very small. However, the division by a very
small number generates bias while our approach is more stable.

\subsection{Dual SDDP for an hydro-thermal generation problem}

We repeat the experiments of Section \ref{sec:inventorysim} for the
Brazilian interconnected power system problem discussed  in \cite{pythonlibrarysddp} for $T=12$ stages
and $N_t=50$ inflow realizations for every stage. These realizations are
obtained calibrating log-normal distributions for each month of the year
using historical data of inflows and sampling from these distributions.
The data used for these simulations (including the inflow scenarios) is available on Github\footnote{{\tt{{\url{https://github.com/vguigues/Primal_SDDP_Library_Matlab}}}}}.

We solve this problem using Primal SDDP and Dual SDDP with penalizations. For this variant of Dual SDDP, a general
procedure to define sequences of penalizations
$(v_{t k})$ ensuring  convergence of the corresponding
Dual SDDP method is to take
$v_{t k}=\gamma_0 \alpha^{k-1}{\textbf{e}}$, $k \geq 1$, $t=2,\ldots,T$, with
$\alpha>1$, $\gamma_0>0$. For numerical reasons, we also take a large upper bound $U$ for these sequences and use
\begin{equation}\label{pensequence}
v_{t k} =  \min(U,\gamma_0 \alpha^{k-1}){\textbf{e}},\;k \geq 1, t=2,\ldots,T.
\end{equation}
We consider three variants of Dual SDDP: for the first variant, denoted by {\tt{Dual SDDP 1}},
$v_{t k}$  are as in \eqref{pensequence} with $\gamma_0=10^4$, $\alpha=1.3$, $U=10^{10}$.
To illustrate the fact that for constant sequences $v_{t k}=\gamma_0$, Dual SDDP converges (resp. does not converge)
for sufficiently large constants $\gamma_0$ (resp. sufficiently small constants $\gamma_0$)
we also define two other variants corresponding to  $U=+\infty$, $\gamma_0=10^9$, $\alpha=1$,
and $U=+\infty$, $\gamma_0=10^6$, $\alpha=1$, in \eqref{pensequence}, respectively denoted by
{\tt{Dual SDDP 2}} and {\tt{Dual SDDP 3}}.

\begin{figure}
 \centering
 \begin{tabular}{cc}
\includegraphics[scale=0.49]{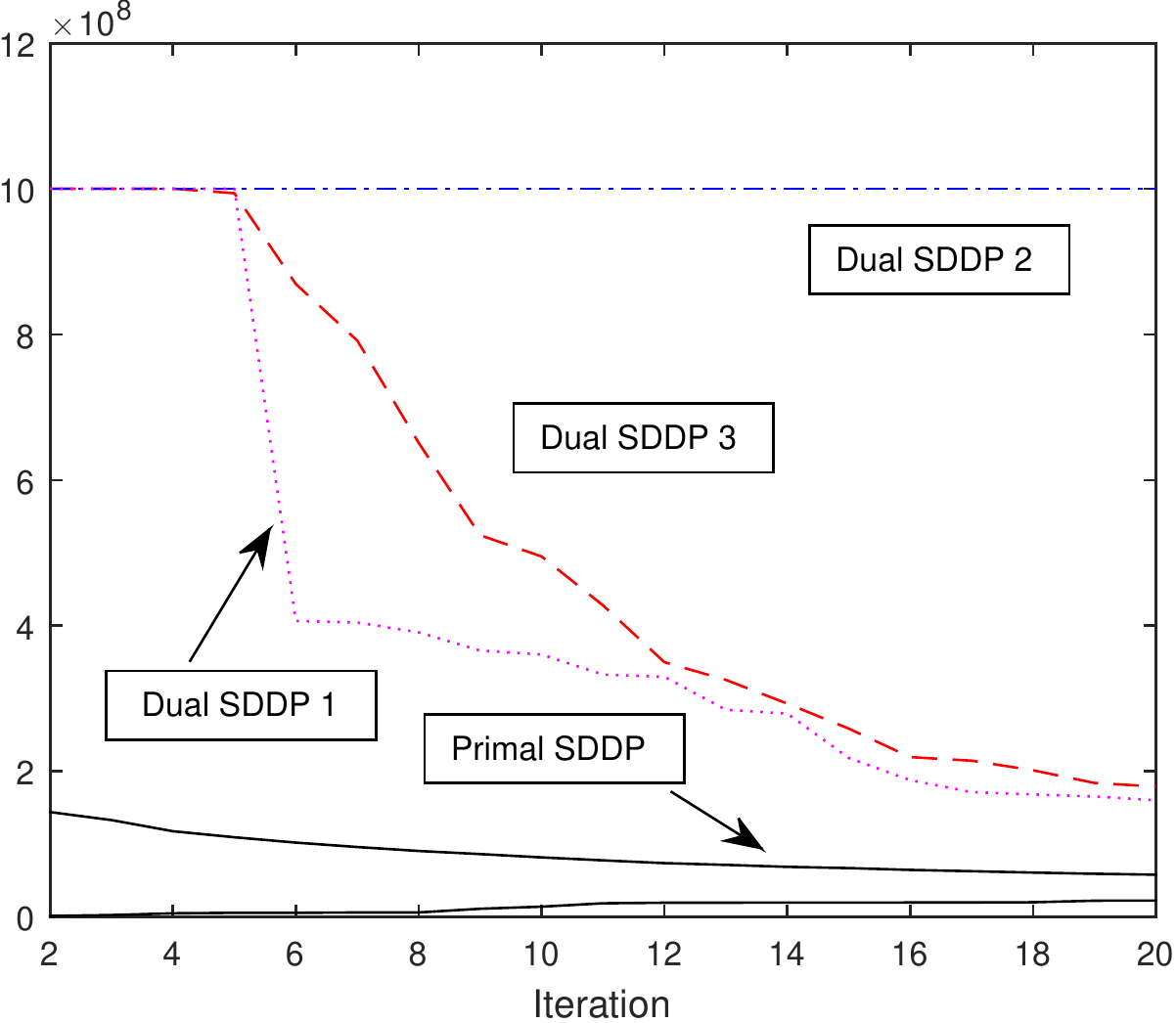}
 &
\includegraphics[scale=0.49]{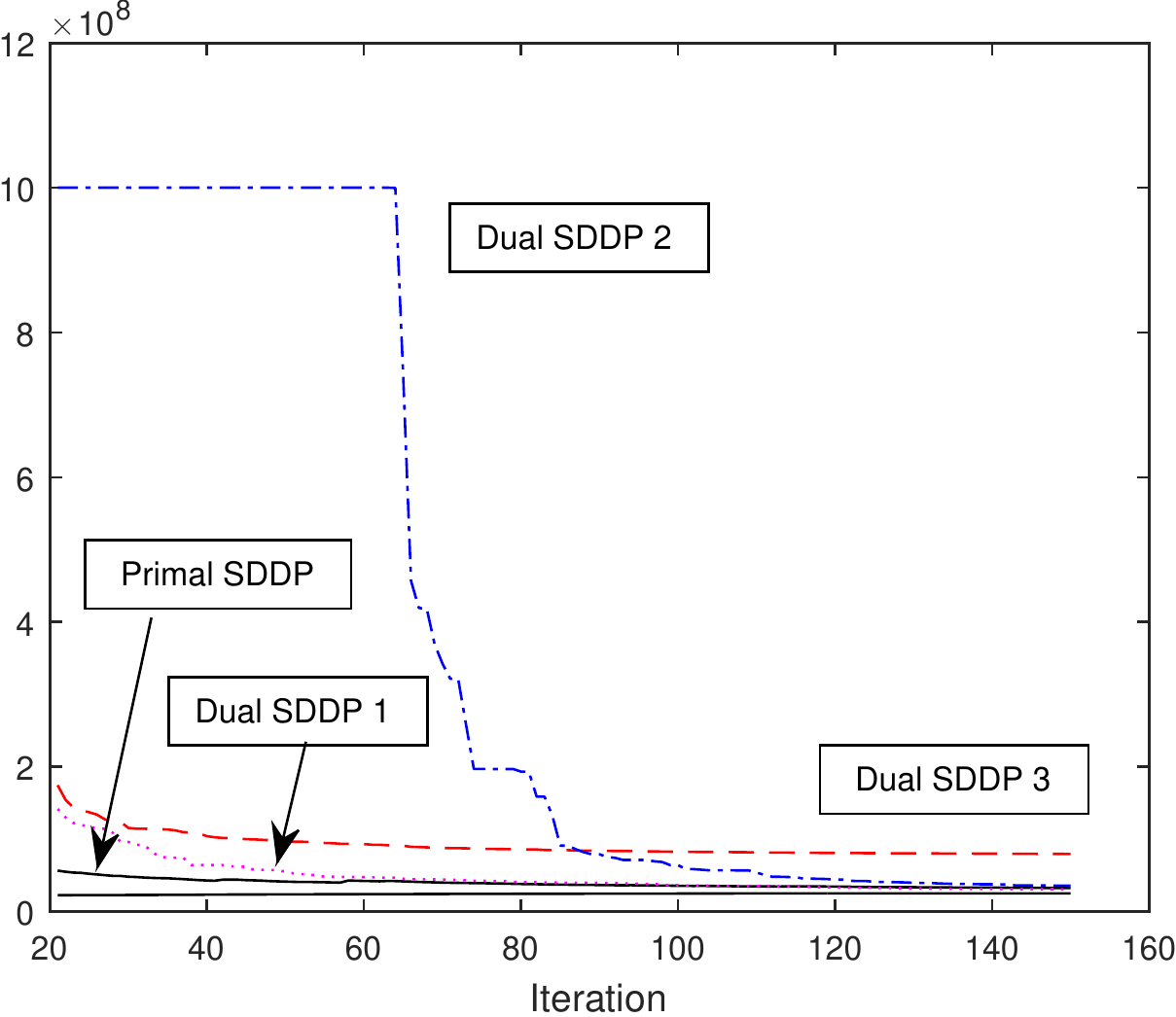}\\
 \end{tabular}
\begin{tabular}{c}
\includegraphics[scale=0.5]{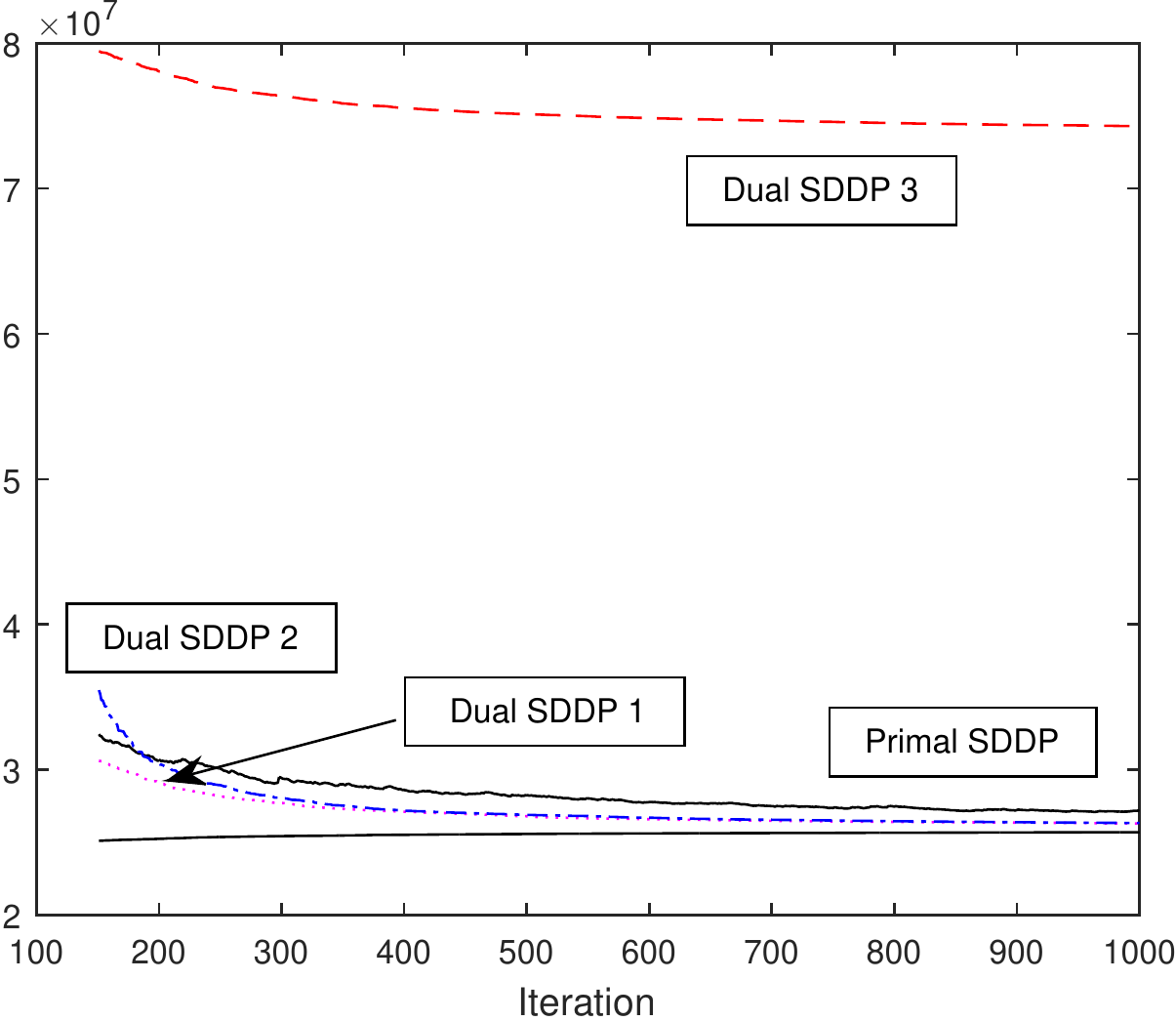}
 \end{tabular}
\caption{Top left: upper and lower bounds computed by Primal SDDP
and upper bounds computed by {\tt{Dual SDDP 1}}, {\tt{Dual SDDP 2}}, and {\tt{Dual SDDP 3}},
for the first 20 iterations for an instance of the hydro-thermal problem with $T=12$, $N_t=50$.
Top right: same outputs for iterations 21,$\ldots,150$.
Bottom: same outputs for iterations 151,$\ldots,1000$. }\label{figboundsht}
\end{figure}

We run Dual SDDP for 1000 iterations and Primal SDDP for 3000 iterations. The evolution of the upper and
lower bounds computed by the methods for the first 1000 iterations is given in Figure \ref{figboundsht}.\footnote{The upper bounds for Primal SDDP are computed as explained in Section \ref{sec:inventorysim}.}

More precisely, the values of these bounds for iterations 2, 5, 10, 50, 100, 150, 200, 250, 300, 350, 400, 1000,
and 3000 are reported in Table \ref{tablebounditerht}. We observe that parameter $\gamma_0$ for {\tt{Dual SDDP 3}}
is too small to allow this method to converge to the optimal value of the problem whereas the other two variants
{\tt{Dual SDDP 1}} and {\tt{Dual SDDP 2}} of Dual SDDP converge. Naturally, these methods start with
large upper bounds but after a few tens of iterations the upper bounds with {\tt{Dual SDDP 1}} and {\tt{Dual SDDP 2}}
are better than the upper bound computed by Primal SDDP. In particular, it is interesting to notice
that the best (lowest) upper bounds are obtained with the variant of Dual SDDP that uses adaptive
penalizations, i.e., penalizations that increase with the number of iterations before reaching value $U$
in \eqref{pensequence}.

\begin{table}
\centering
\begin{tabular}{|c|c|c|c|c|c|}
\hline
Iteration & \begin{tabular}{c}{\tt{Primal}}\\{\tt{SDDP}}\\Lb \end{tabular}&  \begin{tabular}{c}{\tt{Primal}}\\{\tt{SDDP}}\\Ub\end{tabular}  & {\tt{Dual SDDP 1}} & {\tt{Dual SDDP 2}} & {\tt{Dual SDDP 3}}\\
\hline
2  & $1.317$ & $143.98$  & $1000.2$ & $1000.2$&$1000.2$\\
\hline
5  & $5.5588$ & $109.36$  & $1000.2$ &  $1000.2$&$994.04$\\
\hline
10 & $14.032$ &     $81.728$  & $360.40$   &  $1000.2$& $495.08$ \\
\hline
50     & $23.670$ &  $41.346$  & $54.999$  & $1000.2$& $96.720$ \\
\hline
100   & $24.787$ &  $35.502$  & $36.322$  & $64.072$ & $82.494$\\
\hline
150   & $25.111$ &  $32.447$  & $30.685$  & $35.595$& $79.465$\\
\hline
200   & $25.249$ &  $30.672$  & $29.076$  & $30.404$&$78.059$ \\
\hline
250   & $25.374$ &  $30.079$  & $28.215$  & $28.943$& $76.917$\\
\hline
300   & $25.436$ &  $29.434$  & $27.710$  & $28.030$& $76.344$\\
\hline
350   & $25.477$ &  $29.014$  & $27.309$  & $27.532$& $75.852$\\
\hline
400   & $25.526$ &  $28.626$  & $27.110$  & $27.188$& $75.526$\\
\hline
1000  & $25.703$ &  $27.175$  & $26.304$ & $26.335$  &$74.292$\\
\hline
3000  & $25.798$ &  $26.883$  & -  &  - &\\
\hline
\end{tabular}
\caption{For an instance of the hydro-thermal problem with $T=12$, $N_t=50$,
lower bound Lb and upper bound Ub computed by Primal SDDP and upper bounds computed
by variants of Dual SDDP along iterations. All costs have been divided
by $10^6$ to improve readability.}\label{tablebounditerht}
\end{table}

We also report in Table \ref{table1htiter} the relative error
$\frac{   {\tt{Upper}}_M(i)-{\tt{Lower}}_{{\tt{SDDP}}}(i)}{{\tt{Upper}}_M(i)}$ for iterations $i=100$, $200$, $300$, $400$, $500$, $800$, and $1000$
for all methods M where ${\tt{Upper}}_M(i)$ and ${\tt{Lower}}_{{\tt{SDDP}}}(i)$
are respectively the upper bound computed by method M at iteration $i$
and the lower bound computed by Primal SDDP at iteration $i$. For iterations 300 on, the relative error
is much smaller with variants of Dual SDDP, meaning that Primal SDDP overestimates the optimality gap.

\begin{table}
\centering
\begin{tabular}{|c|c|c|c|}
\hline
Iteration & {\tt{Primal SDDP}} & {\tt{Dual SDDP 1}} & {\tt{Dual SDDP 2}} \\
\hline
100  & 0.30 & 0.32 & 0.61\\
\hline
200  &  0.18& 0.13 &0.17\\
\hline
300  & 0.14  & 0.08 &0.09\\
\hline
400  &  0.11 & 0.06&0.06\\
\hline
500 &  0.09  & 0.05 &0.05\\
\hline
800 & 0.07  & 0.03 &0.03\\
\hline
1000 & 0.05   & 0.02 &0.02\\
\hline
\end{tabular}\\
\caption{Relative error as a function of the number of iterations for
{\tt{Primal SDDP}}, {\tt{Dual SDDP 1}}, and {\tt{Dual SDDP 2}}.
}\label{table1htiter}
\end{table}

\begin{figure}
 \centering
\includegraphics[scale=0.65]{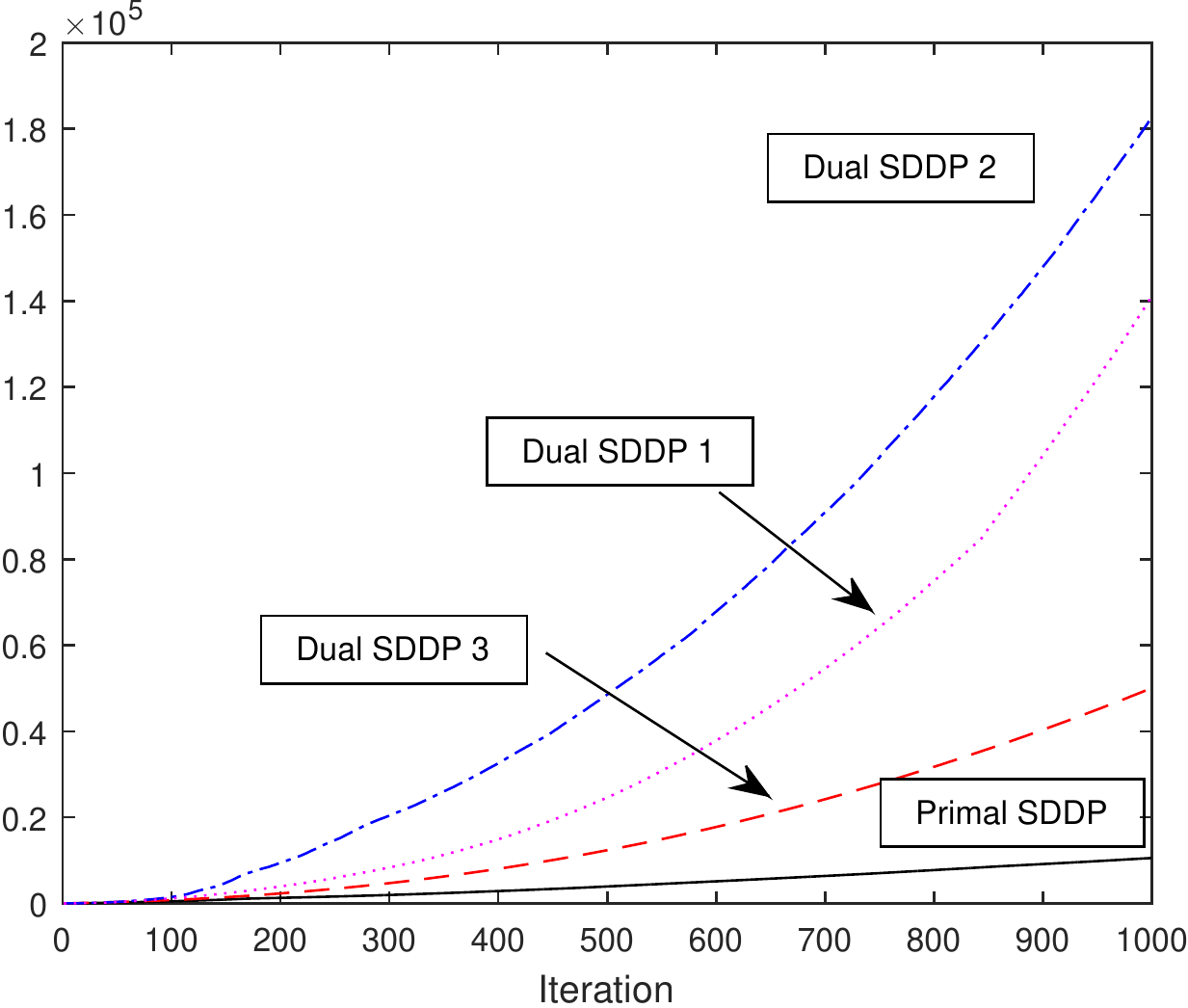}
\caption{Cumulative CPU time for Primal SDDP, Dual SDDP 1, Dual SDDP 2, and Dual SDDP 3.}\label{fighttime}
\end{figure}

However, each iteration of Dual SDDP takes more time as can be seen in
Figure \ref{fighttime} which reports the cumulative CPU time for all methods.
More precisely, running Dual and Primal SDDP in parallel, we can compute the time needed
to obtain a solution of relative accuracy $\varepsilon$ using the standard stopping criterion
for Primal SDDP (see \cite{shapsddp}) or using the lower bound from Primal SDDP and the upper bound
from Dual SDDP, and computing the relative error obtained with these bounds each time a new bound
(either lower bound or upper bound) is computed. The results are reported in Table \ref{table1htprecision}.
We see that due to the fact that Dual SDDP iterations are more time consuming, for all relative accuracies but one,
the use of the stopping criterion based on Dual SDDP upper bounds requires more computational bulk.
From this experiment, performed on a larger problem (in terms
of size of the state vector and number of control variables for each stage) than the inventory problem
of Section \ref{sec:inventorysim}, it seems that the use of Dual SDDP
for a stopping criterion of Primal SDDP will decrease the overall computational
bulk only for small problems (having a limited to small number of controls, state variables, and
scenarios).

\begin{table}
\centering
\begin{tabular}{|c|c|c|c|}
\hline
$\varepsilon$ & {\tt{Primal SDDP}} & {\tt{Dual SDDP 1}} & {\tt{Dual SDDP 2}} \\
\hline
0.3  & 515 & 1 042 &4 133\\
\hline
0.2  & 1 167 & 1 895 &7 446\\
\hline
0.15  & 1 659  & 2 910 &9 882\\
\hline
0.1  &  3 168 &5 114 &16 387\\
\hline
0.075 &  5 359  & 8 003 &22 457\\
\hline
0.05 &  11 124  & 15 738 &35 113\\
\hline
0.04 &  45 391  & 23 449  &51 381\\
\hline
\end{tabular}\\
\caption{Time (in seconds) needed to obtain a solution of relative accuracy $\varepsilon$
with Primal SDDP and variants of Dual SDDP for an instance of the hydro-thermal problem.}\label{table1htprecision}
\end{table}

\if{
Finally, similarly to Figure \ref{fig3error},
as an evidence of the fact that RCR still does not hold for the dual of the hydro-thermal problem, we
observed that the maximal
 and mean values
 of  $\|\zeta_t^{k}\|_1$
along iterations, where $\zeta_t^{k}$ is an optimal
value
of $\zeta_t$ in \eqref{forwarddualsddppen} for iteration $k$, are positive for some stages.
}\fi

\addcontentsline{toc}{section}{References}
\bibliographystyle{plain}
\bibliography{Biblio}

\setcounter{equation}{0}
\section{Appendix}
\label{appen}
In this Appendix, we prove
Lemma \ref{boundedpi}, Proposition \ref{proofubfeassddp}, and Theorem \ref{thconvdualsddppenal}.

We first need more notation.
We introduce the sequence of  functions
\begin{equation}\label{defvhatTbis}
{\overline V}_T^k (\pi_{T-1}):=
\begin{array}{l}
\max\limits_{\pi_{T 1},\ldots,\pi_{T N_T}, \zeta_T} \; \displaystyle \sum\limits_{j=1}^{N_T} p_{T j}
 b_{T j}^\top \pi_{T j}-v_{T k}^\top \zeta_T\\
A_{T}^{\top} \pi_{T j} \leq  c_{T},\;j=1,\ldots,N_T,\\
A_{T-1}^{\top} \pi_{T-1} + \displaystyle \sum_{j=1}^{N_T}
p_{T j} B_{T j}^{\top} \pi_{T j} \leq c_{T-1}
+ \zeta_T,\\
\zeta_T \geq 0, {\underline \pi}_T \leq \pi_{T j} \leq {\overline \pi}_{T},\;j=1,\ldots,N_T,
\end{array}
\end{equation}
and for $t=2,\ldots,T-1$, the sequence of functions
\begin{equation}\label{defvhattbis}
{\overline V}_t^k(\pi_{t-1}):=
\begin{array}{l}
\max\limits_{\pi_{t 1},\ldots,\pi_{t N_t}, \zeta_t} \; \displaystyle \sum\limits_{j=1}^{N_t} p_{t j} \left(
  b_{t j}^\top \pi_{t j} + {V}_{t+1}^k (\pi_{t j}) \right) - v_{t k}^\top \zeta_t\\
A_{t-1}^{\top} \pi_{t-1}+ \displaystyle \sum_{j=1}^{N_t} p_{t j} B_{t j}^{\top} \pi_{t j}  \leq  c_{t-1} + \zeta_t,\\
\zeta_{t} \geq 0,{\underline \pi}_t \leq \pi_{t j} \leq {\overline \pi}_{t},\;j=1,...,N_t.\\
\end{array}
\end{equation}
Due to Assumption (A1) we can
represent the scenarios
for $\xi_1,\xi_2,\ldots,\xi_T$,
by a scenario tree
of depth $T+1$
where the root node $n_0$ associated to a stage $0$ (with decision $x_0$ taken at that
node) has one child node $n_1$
associated to the first stage. We denote by 
$\mathcal{N}$ the set of nodes and for a node $n$ 
of the tree, by $(x_n,\pi_n)$ a primal-dual pair at that 
node and by
$\xi_n$ the realization of process $(\xi_t)$ at node $n$ (this realization $\xi_n$ contains in particular the realizations
$c_n$ of $c_t$, $b_n$ of $b_t$, $A_{n}$ of $A_{t}$, and $B_{n}$ of $B_{t}$).\\

\par {\textbf{Proof of Lemma \ref{boundedpi}.}}
Let $1 \leq t \leq T$ and
let us fix a node $m$
of stage $t$.
Let
$\overline A_m$
such that
constraints
$A_m x_m + B_m x_{F(m)}=b_m$
are rewritten in compact form
$\overline A_m x = b_m$
in terms of vector
$x=(x_n)_{n \in \mathcal{N}}$
of decision variables
in the scenario tree.
The dual function obtained dualizing
the coupling constraints of node
$m$ is given by
$$
\theta(\pi_m) =
\begin{array}{l}
\min \; \mathbb{E}[c^\top x] + \pi_m^\top
(A_m x_m + B_m x_{F(m)}-b_m)\\
x \in \mathcal{S}_m,
\end{array}
$$
for
$
\mathcal{S}_m
=\{
x=(x_n)_{n \in \mathcal{N}} :
x \geq 0 \} \cap \mathcal{A}_m
$
where
$\mathcal{A}_m=\{
x=(x_n)_{n \in \mathcal{N}} :
 A_n x_n +B_n x_{F(n)} = b_n, \forall
n \neq m, n \in \mathcal{N} \}.
$

By Linear Programming Duality, the optimal
value $\mathcal{Q}_1(x_0)$ of primal
problem \eqref{eq-1}
is the optimal value of the dual
problem
\begin{equation}\label{dualpbinit}
\max \{\theta(\pi_m) : \pi_m \in \mathbb{R}^{m_t}\}
\end{equation}
which can clearly be written as
\begin{equation}\label{reformdual}
\mathcal{Q}_1(x_0) = \max_{\pi_m}\{ \theta(\pi_m):
\pi_m = {\overline A}_m x -
{b}_m, x \in \mbox{Aff}(\mathcal{S}_m)\},
\end{equation}
where $\mbox{Aff}(\mathcal{S}_m)$
is the affine hull of  $\mathcal{S}_m$.
We now bound the optimal solutions
of dual problem \eqref{reformdual}.
Since \eqref{dualpbinit}
and \eqref{reformdual} have the same
optimal values, adding these bounds as
constraints on $\pi_m$ in
\eqref{dualpbinit} does not change
its optimal value. Since $\hat x>0$ there is $r>0$ such that
\begin{equation}\label{xhatrelint}
\mathbb{B}(\hat x,r)
\subseteq \{x \geq 0\}.
\end{equation}
We argue that
$\mbox{Aff}(\mathcal{S}_m)=\mathcal{A}_m$.
Indeed, the inclusion
$\mbox{Aff}(\mathcal{S}_m) \subseteq \mathcal{A}_m$ is clear.
Now if $x \in \mathcal{A}_m$ then if $x = \hat x$
we have that $x \in \mathcal{S}_m \subseteq
\mbox{Aff}(\mathcal{S}_m)$ and if $x \neq \hat x$,
recalling
that $\hat x \in \mathcal{A}_m$ satisfies
\eqref{xhatrelint} we have that
$$
y:=\hat  x + \frac{r}{2} \frac{x-\hat x}{\|x-\hat x\|}
\in \mathcal{A}_m \cap \mathbb{B}(\hat x,r)
\subseteq \mathcal{S}_m.
$$
Therefore $x$ belongs to the line that contains $y$ and $\hat x$ with
$y, \hat x$ belonging to
$\mathcal{S}_m$ which implies $x \in \mbox{Aff}(\mathcal{S}_m)$
and
$\mbox{Aff}(\mathcal{S}_m)=\mathcal{A}_m$.

It follows that
$$
\mathbb{B}(\hat x,r)
\cap \mbox{Aff}(\mathcal{S}_m)
= \mathbb{B}(\hat x,r)
\cap \mathcal{A}_m \subseteq \mathcal{S}_m
$$
and that
 there is $\rho_*(m)>0$
such that
$$
\mathbb{B}(0,\rho_*)
\cap (\overline{A}_m \mathcal{A}_m - b_m)
\subseteq \overline A_m
(\mathbb{B}(\hat x,r) \cap \mathcal{A}_m)
-b_m.
$$
Let $\bar \pi_m$ be
an optimal solution of problem \eqref{reformdual}
and let $z=0$ if $\bar \pi_m=0$
and
$z=-\frac{\bar \pi_m}{\|\bar \pi_m\|_2}\rho_*$
otherwise.
Observe that
$z \in \mathbb{B}(0,\rho_*)
\cap (\overline{A}_m \mathcal{A}_m - b_m)
$
and therefore
$z \in \overline A_m
(\mathbb{B}(\hat x,r) \cap \mathcal{A}_m)
-b_m \subseteq
\overline A_m \mathcal{S}_m-b_m$
and $z$ can be written
$z={\overline A}_m \tilde x - b_m$
for $\tilde x \in \mathbb{B}(\hat x,r) \cap \mathcal{S}_m$. It follows
that
$$
\begin{array}{lcl}
\mathcal{Q}_1(x_0)
=\theta(\bar \pi_m)
& \leq &
\mathbb{E}[c^\top \tilde x] +
\bar \pi_m^{\top}({\overline A}_m \tilde x - b_m)\\
& \leq & \mathbb{E}[c^\top \hat x] + r
\sum_{t=1}^T \mathbb{E}[\| c_t \|_2] +  \bar \pi_m^{\top} z\\
& = &  \mathbb{E}[c^\top \hat x] +
r \sum_{t=1}^T \mathbb{E}[\| c_t \|_2]  - \rho_*(m) \|\bar \pi_m\|_2
\end{array}
$$
which gives for every node
$n$ of stage $t$ that
$$
\|\bar \pi_n\|_2
\leq
\max_{m \in {\tt{Nodes}}(t)}
\frac{
\mathbb{E}[c^\top \hat x]
-\mathcal{Q}_1(x_0)+r \sum_{t=1}^T \mathbb{E}[\| c_t \|_2]}{\rho_*(m)}
$$
with corresponding box constraints
$\underline \pi_t, \overline \pi_t$
where ${\tt{Nodes}}(t)$ are the nodes
of stage $t$. $\hfill \square$\\

\par {\textbf{Proof of Proposition \ref{proofubfeassddp}.}} We show by induction on $k$ that $V_t \leq {V}_t^k$ for $t=2,\ldots,T$. For $k=0$ these relations hold by definition. Assume that for some $k \geq 1$ we have $V_t \leq {V}_t^{k-1}$ for $t=2,\ldots,T$. We show by backward induction on $t$ that $V_t \leq {V}_t^{k}$ for $t=2,\ldots,T$. Observe that for any $\pi_{T-1}$, optimization problem \eqref{defvhatTbis} with optimal value ${\overline V}_T^k(\pi_{T-1})$ is feasible. Indeed, since primal problem \eqref{eq-1} is feasible and has a finite optimal value, the corresponding dual problem is feasible which implies that there is $\pi_{T 1},\ldots,\pi_{T N_T}$ satisfying $A_{T}^{\top} \pi_{T j} \leq  c_{T}$, ${\underline \pi}_T \leq \pi_{T j} \leq {\overline \pi}_{T}$, $j=1,\ldots,N_T$, and for every such points we can find $\zeta_T \geq 0$ satisfying the remaining constraints in \eqref{defvhatTbis}. Therefore ${\overline V}_T^k(\pi_{T-1})$ is finite for every $\pi_{T-1}$ and is the optimal value of the corresponding dual optimization problem, i.e., for any $\pi_{T-1}$ we get $$ {\overline V}_T^k(\pi_{T-1}) = \begin{array}{l} \displaystyle \min_{\alpha, \delta, {\overline{\Psi}}, {\underline{\Psi}}} \; \delta^{\top} (c_{T-1} - A_{T-1}^{\top} \pi_{T-1})  + c_T^{\top}  \sum_{j=1}^{N_T} \alpha_j + \sum_{j=1}^{N_T} {\overline{\Psi}}_j^\top {\overline \pi}_T - \sum_{j=1}^{N_T} {\underline{\Psi}}_j^\top {\underline \pi}_T \\ A_T \alpha_j + p_{T j} B_{T j} \delta - {\underline{\Psi}}_j +{\overline{\Psi}}_j = p_{T j} b_{T j},\;j=1,\ldots,N_T,\\ 0 \leq \delta \leq v_{T k}, \alpha_j, {\underline{\Psi}}_j, {\overline{\Psi}}_j \geq 0,\;j=1,\ldots,N_T. \end{array} $$ Using this dual representation and the definition of $\overline{\theta}_T^k, {\overline{\beta}}_T^k$, we get for every $\pi_{T-1}$: \begin{equation}\label{cutaboveT} \overline{\theta}_T^k + \langle {\overline{\beta}}_T^k , \pi_{T-1} \rangle  \geq {\overline V}_T^k(\pi_{T-1}). \end{equation} Recalling representation \eqref{defvhatTbis} for ${\overline V}_T^k(\pi_{T-1})$, observe that for every $\pi_{T-1} \in \mbox{dom}(V_T)$ we have ${\overline V}_T^k(\pi_{T-1}) \geq V_T(\pi_{T-1})$ whereas for $\pi_{T-1} \notin \mbox{dom}(V_T)$ we have $V_T(\pi_{T-1})=-\infty$ while ${\overline V}_T^k(\pi_{T-1})$ is finite, which shows that for every $\pi_{T-1}$ we have $ {\overline V}_T^k(\pi_{T-1}) \geq V_T(\pi_{T-1})$, which, combined with \eqref{cutaboveT} and the induction hypothesis, gives $$V_T^k(\pi_{T-1}) \geq V_T(\pi_{T-1})$$ for every$\pi_{T-1}$.

Now assume that $V_{t+1}^k(\pi_{t}) \geq V_{t+1}(\pi_{t})$ for all $\pi_t$ for some  $t \in \{2,\ldots,T-1\}$. We want to show that $V_{t}^k(\pi_{t-1}) \geq V_{t}(\pi_{t-1})$ for all $\pi_{t-1}$. First observe that for every $\pi_{t-1}$, linear program \eqref{defvhattbis} with optimal value ${\overline V}_t^k(\pi_{t-1})$ is feasible and has a finite optimal value. Therefore we can express ${\overline V}_t^k(\pi_{t-1})$ as the optimal value of the corresponding dual problem given by \begin{equation}\label{dualpbvtk1} \begin{array}{l} \displaystyle \min_{\delta,\nu,{\overline \Psi}, {\underline \Psi}} \; \delta^{\top}\Big[ c_{t-1} - A_{t-1}^{\top} \pi_{t-1}\Big] + \sum_{i=0}^{k} {\overline{\theta}}_{t+1}^i \sum_{j=1}^{N_t} \nu_i(j) + \sum_{j=1}^{N_t} {\overline \Psi}_j^\top {\overline \pi}_t -  \sum_{j=1}^{N_t} {\underline \Psi}_j^\top {\underline \pi}_t\\ \displaystyle p_{t j} B_{t j} \delta - \sum_{i=0}^k \nu_i(j)  {\overline{\beta}}_{t+1}^i - {\underline \Psi}_j + {\overline \Psi}_j = p_{t j} b_{t j},\;j=1,\ldots,N_t,\\ \displaystyle \sum_{i=0}^k \nu_i( j ) = p_{t j}, {\underline \Psi}_j, {\overline \Psi}_j \geq 0,\;j=1,\ldots,N_t,\\ \nu_0,\ldots,\nu_k \geq 0,0 \leq \delta \leq v_{t k}. \end{array} \end{equation} Using this representation of $\overline V_t^k$ and the definition of $\overline{\theta}_t^k, {\overline{\beta}}_t^k$, we obtain for every $\pi_{t-1}$: \begin{equation}\label{cutabovet} \overline{\theta}_t^k +  \langle {\overline{\beta}}_t^k , \pi_{t-1} \rangle \geq {\overline V}_t^k(\pi_{t-1}). \end{equation} Next, recalling representation \eqref{defvhattbis} for ${\overline V}_t^k(\pi_{t-1})$ and the induction hypothesis, we get \begin{equation}\label{induction2} {\overline V}_t^k(\pi_{t-1}) \geq {\widehat V}_t^k(\pi_{t-1}) \end{equation} where $$ {\widehat V}_t^k(\pi_{t-1}):= \begin{array}{l} \max\limits_{\pi_{t 1},\ldots,\pi_{t N_t}, \zeta_t} \; \displaystyle \sum\limits_{j=1}^{N_t} p_{t j} \left( b_{t j}^\top \pi_{t j} + {V}_{t+1}(\pi_{t j}) \right) - v_{t k}^\top \zeta_t\\ A_{t-1}^{\top} \pi_{t-1}+ \displaystyle \sum_{j=1}^{N_t} p_{t j} B_{t j}^{\top} \pi_{t j}  \leq  c_{t-1} + \zeta_t,\\ \zeta_{t} \geq 0,{\underline \pi}_t \leq \pi_{t j} \leq {\overline \pi}_{t},\;j=1,...,N_t.\\ \end{array} $$ Similarly to the induction step $t=T$, for every $\pi_{t-1}$, we have \begin{equation}\label{induction3} {\widehat V}_t^k(\pi_{t-1}) \geq V_t(\pi_{t-1}). \end{equation} Combining \eqref{cutabovet}, \eqref{induction2}, and \eqref{induction3} with the induction hypothesis, we obtain $V_{t}^k(\pi_{t-1}) \geq V_{t}(\pi_{t-1})$ for all $\pi_{t-1}$ which achieves the proof of the induction step $t$.

In particular $V_2^{k-1} \geq V_2$ which implies that $V^{k-1}$ is greater than or equal to the optimal value of dual problem \eqref{eq-2} which is also, by linear programming duality, the optimal value of primal problem \eqref{eq-1}.$\hfill \square$\\

The proof of Theorem \ref{thconvdualsddppenal} is based on the following lemma:

\begin{lemma}
\label{lem-conver}
Suppose that the multistage  problem  \eqref{eq-1} has a finite optimal value. Then for sufficiently large values of the components of vectors $v_t$,  in the dynamic equations \eqref{eq-r2}, the optimal   value of the multistage problem defined by these dynamic equations coincides with the optimal value of the original problem \eqref{eq-1}.
\end{lemma}

 {\bf Proof.} 
  As it was already mentioned,  since it is assumed that the number of scenarios is finite,
 we can view problem \eqref{eq-1} as a large  linear program (deterministic equivalent) written under the form 
 \begin{equation}\label{linprimalpr-1}
\min_{x} c^\top x\;\;{\rm s.t.}\; 
\A x =b,\;x \geq 0.  
 \end{equation}
Also since \eqref{eq-1} has a finite optimal
 value, it has a nonempty set of optimal solutions
 and there is a bounded optimal solution 
of \eqref{linprimalpr-1}.
Let us fix such an optimal solution $\bar x$. 
We have that problem \eqref{linprimalpr-1}  can
 be written
 \begin{equation}\label{linprimalpr-2}
\min_{x} c^\top x\;\;{\rm s.t.}\; 
\A x =b,\;0 \leq x \leq \bar x.  
 \end{equation}
 The dynamic programming equations \eqref{eq-3a-1} - \eqref{eq-7-1} represent the standard dual of \eqref{eq-1}. We can also  think about that dual as a large linear programming problem of the form
 (this is the dual of \eqref{linprimalpr-1}):
 \begin{equation}\label{lindualpr-1}
\max_{\pi} b^\top \pi\;\;{\rm s.t.}\; 
\A^\top \pi \le c.  
 \end{equation}
Similarly the deterministic equivalent of  penalized dynamic equations \eqref{eq-r2} can be written as:
 \begin{equation}\label{lindualpr-2}
\max_{\pi,\zeta} b^\top  \pi-v^\top \zeta\;\;{\rm s.t.}\;
\A^\top \pi \le c+\zeta, \zeta \geq 0.
 \end{equation}
Next, from optimality conditions of linear programs,
$(x,\pi)$ is an optimal primal-dual pair
for \eqref{linprimalpr-1}-\eqref{lindualpr-1}
if and only if 
 \begin{equation}\label{lindualpr-3}
x^\top(\A^\top \pi -c)=0,\;
\A x =b,\;x \geq 0,\; 
\A^\top \pi \leq c.
 \end{equation}
 The corresponding optimality conditions for \eqref{lindualpr-2} are
 \begin{equation}\label{lindualpr-4}
 x^\top(\A^\top \pi-c-\zeta)-\zeta^\top \gamma=0,\;
\A^\top \pi \leq c+ \zeta,\;
\zeta \geq 0,\;
\A x = b,\;x \geq 0,\;\gamma \geq 0,\;
x=v-\gamma. 
\end{equation}
Now let $\bar{\pi}$ be an optimal dual solution, i.e.,
an optimal solution of \eqref{lindualpr-1}.
Then \eqref{lindualpr-3}
is satisfied with 
$(x,\pi)=(\bar x, \bar \pi)$.
It follows  that if $v\ge \bar x$, then 
$(x,\pi,\zeta,\gamma)=(\bar x,\bar \pi,0,v-\bar x)$  with  
$\zeta=0$ satisfies \eqref{lindualpr-4}, and hence 
$(\bar \pi,\bar \zeta)=(\bar \pi,0)$
is an optimal solution of  
\eqref{lindualpr-2} showing that the optimal value
of \eqref{lindualpr-2} is
$b^\top \bar \pi=c^\top \bar x$, i.e., the
optimal value of \eqref{linprimalpr-1}.
We obtain  that 
for $v\ge \bar x$, the optimal values of problems  \eqref{lindualpr-1} and \eqref{lindualpr-2} do coincide.\footnote{Observe that the dual of \eqref{lindualpr-2}
is given by
$$
\min_{x} c^\top x\;\;{\rm s.t.}\; 
\A x =b,\;0 \leq x \leq v,
$$
and for $v \geq \bar x$, this linear program has the same
optimal value as \eqref{linprimalpr-2}, which, as we have seen, is equivalent to primal problem \eqref{eq-1}.}$\hfill \square$\\

\par {\textbf{Proof of Theorem \ref{thconvdualsddppenal}.}} Dual SDDP with
penalizations is SDDP applied to Dynamic Programming equations  
corresponding to a linear program with finite optimal value, satisfying relatively complete recourse
with discrete uncertainties of finite support. Since samples $\tilde \xi_t^k$
in Dual SDDP with penalizations are independent, 
we can follow the convergence proof of SDDP for linear programs
from \cite{philpot} 
to obtain that  $V^k$ converges
to the optimal value of the penalized linear programs, which,
by Lemma \ref{lem-conver} (observe that the Lemma can be applied
since $\lim_{k \rightarrow +\infty} v_{t k}=+\infty$),
is the optimal value  of \eqref{eq-1}.$\hfill \square$

\if{

Let us take an arbitrary realization of the algorithm and let us show that
$V^k$ converges to
the optimal value of
\eqref{eq-1}
for that realization. Define $V_{T+1}$ and
${\overline{V}}_{T+1}^k$ for $k \geq 0$
on $[\underline{\pi}_T,\overline{\pi}_T]$
to be the null function.
We show that
\begin{equation}\label{inductionkforward}
\lim_{k \rightarrow +\infty} V_t^k(\pi_n^k)-V_t(\pi_n^k)=0
\end{equation}
almost surely
for all $t=2,\ldots,T+1$ and all node $n$ of stage
$t-1$ where $\pi_n^k$ is computed simulating the dual policy obtained
using Dual SDDP with penalizations replacing $V_t$ by $V_t^{k-1}$.
The conclusion of the theorem is an immediate consequence of
\eqref{inductionkforward} written for $t=2$.

We show \eqref{inductionkforward} by backward induction on $t$.
The relation holds for $t=T+1$. Assume it holds for
$t+1$ with $t \in \{2,\ldots,T\}$ and let us show it holds for stage $t$.
Take a node $n$ of stage $t-1$.
It is enough to show that
$\lim_{k \rightarrow +\infty, k \in \mathcal{S}_n} V_t^k(\pi_n^k)-V_t(\pi_n^k)=0$
where $\mathcal{S}_n=\{k: (\xi_1,\tilde \xi_2^k,\ldots,\tilde \xi_{t-1}^k) =
\xi_{[n]}\}$ with $\xi_{[n]}$ the scenario from node $n_1$
to node $n$. Indeed, as in \cite{lecphilgirar12, guiguessiopt2016}, from this relation, \eqref{inductionkforward}
can then be shown by contradiction using the
Strong Law of Large Numbers and the assumption of independence
of $\xi_t^k$. Therefore let $k \in \mathcal{S}_n$.
Since $V_t^{k} \leq V_t^{k-1}$ (by definition
of $V_t^k)$
and $v_{t k}\geq v_{t k-1}$,
we have
\begin{equation}\label{decVtkbar}
{\overline V}_t^k \leq {\overline V}_t^{k-1}.
\end{equation}
Since $k \in \mathcal{S}_n$, we have
$\overline \theta_t^k
+ \langle \overline \beta_t^k , \pi_{n}^k \rangle
= \overline V_t^k(\pi_{n}^k)$ and therefore
for $0 \leq j \leq k$ we get
$
\overline \theta_t^k
+ \langle \overline \beta_t^k , \pi_{n}^k \rangle
= \overline V_t^k(\pi_{n}^k)
 \leq  {\overline V}_t^j(\pi_{n}^k)
\leq \overline \theta_t^j
+ \langle \overline \beta_t^j , \pi_{n}^k \rangle,
$
which implies, by definition of $V_t^k$,
\begin{equation}\label{decVtkbar2}
V_t^k(\pi_{n}^k) = {\overline V}_t^k(\pi_{n}^k).
\end{equation}
Therefore for $k \in \mathcal{S}_n$, we have
\begin{equation}\label{proofconvdsddp}
0  \leq  V_t^k(\pi_n^k)-V_t(\pi_n^k)    \stackrel{\eqref{decVtkbar2}}{=}  {\overline V}_t^k(\pi_n^k)-V_t(\pi_n^k)
 \leq  {\overline V}_t^{F, k-1}(\pi_n^k)-V_t(\pi_n^k)
\end{equation}
where ${\overline V}_t^{F, k-1}(\pi_n^k)$ is the optimal value of
\begin{equation}\label{VFtk}
\begin{array}{l}
\max\limits_{\pi_{m}, \zeta_t} \; \displaystyle \overline F_t^{k-1}((\pi_m)_{m \in C(n)})= \sum\limits_{m \in C(n)} p_{m} \left(
  b_{m}^\top \pi_{m} + {V}_{t+1}^{k-1} (\pi_{m}) \right) - v_{t k}^\top \zeta_t\\
A_{t-1}^{\top} \pi_{n}^k + \displaystyle \sum_{m \in C(n)} p_{m} B_{m}^{\top} \pi_{m}  \leq  c_{t-1} + \zeta_t,\\
\zeta_{t} \geq 0,{\underline \pi}_t \leq \pi_{m} \leq {\overline \pi}_{t},\;m \in C(n).\\
\end{array}
\end{equation}
The problem above is the problem
solved in the forward pass for
iteration $k$ stage $t$.
By first order optimality conditions,
$(\pi^k=(\pi_m^k)_{m \in C(n)},\zeta_t^k)$ is
an optimal primal solution of this problem
and $(\delta_t^k, \lambda_t^k, {\underline \Psi}_m^k, {\overline \Psi}_m^k)$\footnote{We use for dual variables the same notation as in dual problem \eqref{dualpbvtk1}} is an optimal dual
solution if and only they solve the following
KKT system:
\begin{equation}\label{kkt1}
\begin{array}{l}
A_{t-1}^\top \pi_n^k + \sum_{m \in C(n)} p_m
B_n^{\top} \pi_m^k \leq c_{t-1}+\zeta_t^k, \zeta_t^k \geq 0, \underline \pi_t^k \leq \pi_m^k \leq \overline \pi_t,\forall m \in C(n),\\
\delta_t^k, \lambda_t^k, \underline \Psi_m^k, \overline \Psi_m^k \geq 0,\\
0=-p_m(b_m+s_m)+p_mB_m \delta_t^k - \underline \Psi_m^k
+ \overline \Psi_m^k,\\
0=v_{t k} -\delta_t^k - \lambda_t^k,\\
\lambda_t^k(i)\zeta_t^k(i)=0,\underline{\Psi}_m(i)
(\underline{\pi}_t(i)-\pi_m(i))=0,
\overline{\Psi}_m(i) (\pi_m(i)-\overline{\pi}_t(i))=0,\forall i,\\
\delta_t^k(i)(A_{t-1}^\top \pi_n^k + \sum_{m \in C(n)} p_m
B_n^{\top} \pi_m^k - c_{t-1}-\zeta_t^k)(i)=0, \forall i,
\end{array}
\end{equation}
for some $s_m \in \partial V_{t+1}^{k-1}(\pi_m^k)$.
Let
$$\mathcal{X}_t(\pi_n)=
\{(\pi_m):
A_{t-1}^{\top} \pi_{n} + \sum_{m \in C(n)} p_{m} B_{m}^{\top} \pi_{m}  \leq  c_{t-1}, {\underline \pi}_t \leq \pi_{m} \leq {\overline \pi}_{t},\;m \in C(n)\}.$$
Let also $({\hat \pi}_m^k)_{m \in C(n)}$
be an optimal solution of
\begin{equation}\label{VFtkbis}
\begin{array}{l}
\max\limits_{\pi_{m}} \; \displaystyle
F_t^{k-1}((\pi_m)_{m \in C(n)})=
 \sum\limits_{m \in C(n)} p_{m} \left(
  b_{m}^\top \pi_{m} + {V}_{t+1}^{k-1} (\pi_{m}) \right)\\
(\pi_m)_{m \in C(n)} \in  \mathcal{X}_t(\pi_n)
\end{array}
\end{equation}
and let
$(\hat \delta_t^k, {\hat{{\underline \Psi}}}_m^k,
{\hat{{\overline \Psi}}}_m^k)$
be an optimal dual solution (using a notation similar as before
for dual variables).
We partition $\mathcal{S}_n$ in $D_n \cup {\overline{D}}_n$
where $D_n$ is the set of iterations in
$\mathcal{S}_n$ such that
$F_t(\pi_n^k)$ is nonempty.
We first argue that due to
the compactness of
$[\underline{\pi}_t,\overline{\pi}_t]$ and the equicontinuity of
$V_t^k$ there can only be a finite
number of iterations in $\overline D_n$
so that after some  iteration $k_0$
all iterations in $\mathcal{S}_n$
are in $D_n$.
Next, following the proof
of Lemma \ref{boundedpi}, using the fact that
$\underline{\pi}_t<\overline{\pi}_t$,
and that $F_t^{k-1}\leq F_t^{0}$
with $F_t^{0}$ Lipschitz continuous, we get that
there is $\delta_*$  such that for all
$k$ we have $0 \leq \hat \delta_t^k \leq \delta_*$.
From the first order optimality conditions,
$({\hat \pi}_m^k)_{m \in C(n)}, (\hat \delta_t^k, {\hat{{\underline \Psi}}}_m^k,
{\hat{{\overline \Psi}}}_m^k)$
is an optimal primal-dual solution
of \eqref{VFtkbis} if and only if
they solve the following
KKT system:
\begin{equation}\label{kkt2}
\begin{array}{l}
A_{t-1}^\top \pi_n^k + \sum_{m \in C(n)} p_m
B_n^{\top} \hat \pi_m^k \leq c_{t-1}, \underline \pi_t^k \leq \hat \pi_m^k \leq \overline \pi_t,\forall m \in C(n),\\
\hat \delta_t^k, {\hat{\underline \Psi}}_m^k,
{\hat{\overline{\Psi}}}_m^k \geq 0,\\
0=-p_m(b_m+s_m)+p_m B_m {\hat \delta}_t^k - {\hat{\underline \Psi}}_m^k
+ {\hat{\overline \Psi}}_m^k,\\
(\underline{\pi}_t(i)-\hat \pi_m(i))=0,
{\hat{\overline{\Psi}}}_m(i) ({\hat{\pi}}_m(i)-\overline{\pi}_t(i))=0,\forall i,\\
{\hat{\delta}}_t^k(i)(A_{t-1}^\top \pi_n^k + \sum_{m \in C(n)} p_m
B_n^{\top} \hat \pi_m^k - c_{t-1})(i)=0, \forall i,
\end{array}
\end{equation}
for some $s_m \in \partial V_{t+1}^{k-1}(\pi_m^k)$.
Now due to $\lim_{k \rightarrow +\infty} v_{tk}=+\infty$
we can choose $k_0$ such that
for $k \geq k_0$ we have
$v_{t k}>\delta_*$ and therefore
any primal-dual optimal solution $(
({\hat \pi}_m^k)_{m \in C(n)}, (\hat \delta_t^k, {\hat{{\underline \Psi}}}_m^k,
{\hat{{\overline \Psi}}}_m^k))$
of \eqref{VFtkbis} which satisfies
\eqref{kkt2} provides a primal-dual optimal solution
for \eqref{VFtk} given by
$(\pi_m^k)_{m \in C(n)}=({\hat \pi}_m^k)_{m \in C(n)}$,
$\zeta_t^k=0$, $\lambda_t^k=v_{t k}-\delta_t^k$,
$\delta_t^k=\hat \delta_t^k, {{{\underline \Psi}}}_m^k={\hat{{\underline \Psi}}}_m^k,
{{{\overline \Psi}}}_m^k =
{\hat{{\overline \Psi}}}_m^k)$ (because this primal-dual pair satisfies \eqref{kkt1}). Therefore for $k \geq k_0$
the optimal value $\zeta_t^k$ of $\zeta_t$
in \eqref{VFtk} is 0 which gives
${\overline V}_t^{F, k-1}(\pi_n^k)=\sum_{m \in C(n)} p_m\left[b_m^\top \pi_m^k + V_{t+1}^{k-1}(\pi_m^k)\right]$, which, plugged into
gives \eqref{VFtk}, gives
\begin{equation}
\begin{array}{lcl}
V_t^k(\pi_n^k)-V_t(\pi_n^k)& \leq & \displaystyle
\sum_{m \in C(n)} p_m\left[b_m^\top \pi_m^k + V_{t+1}(\pi_m^k)\right]-V_t(\pi_n^k)\\
&& + \displaystyle  \sum_{m \in C(n)} p_m\left[V_{t+1}^{k-1}(\pi_m^k)-V_{t+1}(\pi_m^k)\right] \\
& \leq & \displaystyle \sum_{m \in C(n)} p_m\left[V_{t+1}^{k-1}(\pi_m^k)-V_{t+1}(\pi_m^k)\right],
\end{array}
\end{equation}
where the last inequality comes from the fact that
$\sum_{m \in C(n)} p_m\left[b_m^\top \pi_m^k + V_{t+1}(\pi_m^k)\right] \leq V_t(\pi_n^k)$
by definition of $V_t$ and feasibility of $\pi_m^k$.

By the induction hypothesis,
$\lim_{k\rightarrow +\infty} V_{t+1}^{k}(\pi_m^k)-V_{t+1}(\pi_m^k)=0$ for all $m \in C(n)$
which, using the equicontinuity of $(V_{t+1}^k)_k$ and the
compacity of $[\underline{\pi}_t, \overline \pi_t]$
implies $\lim_{k\rightarrow +\infty} V_{t+1}^{k-1}(\pi_m^k)-V_{t+1}(\pi_m^k)=0$ for all $m \in C(n)$.
Together with \eqref{proofconvdsddp}
this implies
$\lim_{k\rightarrow +\infty, k \in \mathcal{S}_n} V_t^k(\pi_n^k)-V_t(\pi_n^k)=0$
and achieves the proof.
$\hfill \square$}\fi

\end{document}